\documentclass[12pt]{article}
\usepackage{amssymb}
\usepackage{enumerate}

\usepackage{amsfonts}
\usepackage{amsmath,amsthm}

\newcommand{\oper}[2]{\newcommand{#1}{\mathop{\mathrm{#2}}\nolimits} }
\oper{\tr}{tr} \oper{\adj}{adj} \oper{\Div}{div} \oper{\ad}{ad}
\oper{\Ad}{Ad} \oper{\End}{End} \oper{\Hom}{Hom} \oper{\Aut}{Aut}
\oper{\SO}{SO} \oper{\SP}{Sp} \oper{\SU}{SU} \oper{\GL}{GL}
\oper{\Ann}{Ann} \oper{\T}{T} \oper{\U}{U} \oper{\id}{I}
\oper{\ext}{Ext} \oper{\rank}{rank} \oper{\diag}{Diag}
\oper{\re}{Re} \oper{\im}{Im} \oper{\sign}{Sign}

\oper{\DGA}{DGA}

\def\dbar{\overline\partial}

\newcommand{\CC}{\mathbb{C}}

\newcommand{\RR}{\mathbb{R}}

\def\bw{\wedge}

\newcommand{\lie}[1]{\mathfrak{#1}}

\newcommand{\sbr}[2]{[ {#1} \bullet {#2} ]}

\def\tick{\checkmark}

\newcommand{\om}[1]{\omega^{#1}}
\newcommand{\bom}[1]{{\overline\omega}^{#1}}
\newcommand{\te}[1]{\sigma^{#1}}
\newcommand{\bte}[1]{{\overline\sigma}^{#1}}
\newcommand{\ta}[1]{\tau^{#1}}
\newcommand{\et}[1]{\eta^{#1}}

\newcommand{\abs}[1]{\left\lvert #1\right\rvert}

\newenvironment{spmatrix}{\left(\smallmatrix}{\endsmallmatrix\right)}

\newtheorem{theorem}{Theorem}
\newtheorem{corollary}[theorem]{Corollary}
\newtheorem{proposition}[theorem]{Proposition}
\newtheorem{definition}[theorem]{Definition}
\newtheorem{lemma}[theorem]{Lemma}

\newtheorem{remark}{Remark}
\newcommand{\bproof}{\noindent{\it Proof: }}
\newcommand{\eproof}{\  q.~e.~d. \vspace{0.2in}}

\newtheorem*{ack}{Acknowledgments}

\begin{document}

\title{Differential Gerstenhaber Algebras \\
Associated to Nilpotent Algebras}

\author{ R. Cleyton \thanks{Address: Institute f\"ur Mathematik, Humboldt
    Universit\"at zu Berlin, Unter den Linden 6, 10099 Berlin, Germany
    {cleyton@mathematik.hu-berlin.de}. Partially supported by the
    Department of Mathematics, University of California of Riverside,
    the Junior Research Group ``Special Geometries in Mathematical
    Physics'', of the Volkswagen foundation and the SFB~647
    ``Space--Time--Matter'', of the DFG.} \and Y.~S. Poon \thanks{ Address:
    Department of Mathematics, University of California at Riverside,
    Riverside CA 92521, U.S.A. {ypoon@math.ucr.edu}.} }


\maketitle

\abstract{This article provides a complete description of the
  differential Gerstenhaber algebras of all nilpotent complex
  structures on any real six-dimensional nilpotent algebra. As an
  application, we classify all pseudo-K\"ahlerian complex structures
  on six-dimensional nilpotent algebras whose differential
  Gerstenhaber algebra is quasi-isomorphic to that of the symplectic
  structure. In a weak sense of mirror symmetry, this gives a
  classification of pseudo-K\"ahler structures on six-dimensional
  nilpotent algebras whose mirror images are themselves.}

  \

  \noindent{Keywords: Nilpotent Algebra, Gerstenhaber Algebra, Complex Structure, Symplectic Structure, Deformation, Mirror Symmetry.}

  \

  \noindent{AMS Subject Classification: Primary 32G05. Secondary 32G07, 13D10, 16E45, 17B30, 53D45.}

\section{Introduction}
Nilmanifolds, i.e. compact quotients of simply connected nilpotent Lie
groups, are known to be a rich source of exotic geometry. We are
particular interested in pseudo-K\"ahler geometry and its deformation
theory on these spaces. We initially focus on the complex structures,
and will bring symplectic structures in the picture at the end.

It is a general principle that the deformation theories of complex and
symplectic structures are dictated by their associated differential
Gerstenhaber algebras \cite{GM} \cite{Mer2} \cite{Zhou}.  The
associated cohomology theories are Dolbeault's cohomology with
coefficients in the holomorphic tangent bundle, and de Rham's
cohomology respectively. De Rham cohomology of nilmanifolds is known
to be given by invariant differential forms \cite{Nomizu} and there
are several results for Dolbeault cohomology on nilmanifolds pointing
in the same direction \cite{CF} \cite{CFP}. Therefore, in this paper
we focus on invariant objects, i.e. invariant complex structures and
invariant symplectic forms on nilpotent Lie algebras.

Analysis and classification of invariant complex structures and
pseudo-K\"ahler pairs on six-dimensional nilpotent algebras have been
in progress in the past ten years \cite{BDM} \cite{CF} \cite{CFGU}
\cite{Salamon} \cite{Ugarte}. In particular, it is known that a
complex structure can be part of a pseudo-K\"ahler pair, only if it is
nilpotent \cite{CFU}.

After a preliminary presentation on construction of differential
Gerstenhaber algebra for invariant complex and symplectic structures,
we give two key technical results, Proposition \ref{key technical} and
Proposition \ref{second technical}, describing the restrictive nature
of quasi-isomorphisms in our setting.  We recall the definition of
nilpotent complex structure in Section \ref{sec:nil}.  Numerical
invariants for these complex structures are identified, and used to
refine older classifications.  This in particular allows to identify
the real algebra underlying a set of complex structure equations by
evaluation of the invariants. The results of Section~\ref{sec:nil}
including the invariants of complex structure equations and the
associated underlying real algebras are summarized in Table
\ref{tab:1}.

In Section \ref{sec:f}, we analyze the differential Gerstenhaber
algebra $\DGA({\lie g}, J)$ when a nilpotent complex structure $J$ on
a nilpotent Lie algebra $\lie g$ is given. The invariants of the
complex structure equations dictate the structure of $\DGA({\lie g},
J)$.  Relying on the classification provided in Section \ref{sec:nil}
and Table \ref{tab:1}, and in terms of the same set of invariants, we
establish a relation between the Lie algebra structure of $\lie g$ and
that of $\DGA ({\lie g}, J)$. The total output of Section \ref{sec:f}
is provided in Theorem \ref{invariant theorem} and Table \ref{tab:f1}.
These results demonstrate the phenomenon of ``jumping'' of $\DGA(\lie
g, J)$ as $J$ varies through a family of nilpotent complex structures
on some fixed algebra $\lie g$.  Results are given in Theorem
\ref{finding f} and Table \ref{tab:g-f1'}.  With the aid of
Proposition \ref{second technical}, Theorem \ref{invariant theorem} we
also show that each differential Gerstenhaber algebra $\DGA(\lie g,
J)$ is isomorphic to a differential Gerstenhaber algebra $\DGA(\lie h,
O)$ derived from a certain Lie algebra $\lie h$ and linear isomorphism
$O\colon\lie h\to\lie h^*$.  The result is stated in Theorem \ref{iso
  dga}. However, the map $O$ is not necessarily induced by a
contraction with any symplectic form.  A priori, it may not even be
skew-symmetric.

Finally in Section \ref{sec:appl}, we consider the differential
Gerstenhaber algebra $\DGA({\lie h}, \Omega)$ associated to an
invariant symplectic structure $\Omega$ on a nilpotent algebra $\lie
h$. We shall explain in Section \ref{sec:dga} that $\DGA({\lie h},
\Omega)$ is essentially generated by the Lie algebra structure on
${\lie h}$. This elementary observation, along with the results
established in Section \ref{sec:f} and Table \ref{tab:f1}, allows us
to answer the following question: Which six-dimensional nilpotent
algebra $\lie g$ admits a pseudo-K\"ahler structure $(J, \Omega)$ such
that there is a quasi-isomorphism
\[
\DGA ({\lie g}, J) \longrightarrow \DGA({\lie g}, \Omega) \  ?
\]
The construction of DGAs for complex structures and symplectic
structures is well known (e.g. \cite{Zhou}). It is a key ingredient
in homological mirror symmetry.  Extending the concept of mirror
symmetry, Merkulov considers the notion of \it weak \rm mirror
symmetry \cite{Mer2} \cite{Mer}. In this paper, we call a Lie
algebra ${\lie g}$ with a complex structure $J$ and a Lie algebra
${\lie h}$ with a symplectic structure $\Omega$ a ``\it weak \rm
mirror pair'' if there is a quasi-isomorphism between $\DGA ({\lie
g}, J)$ and $\DGA ({\lie h}, \Omega)$. The aforementioned question
stems from a consideration on when ``self mirror'' occurs. For
four-dimensional nilpotent algebras, the answer could be derived
from results in \cite{Poon}. For six-dimensional nilpotent algebras,
our answer is in Theorem \ref{main}.

\section{Differential Gerstenhaber Algebras}\label{sec:dga}

\subsection{Preliminaries}
\begin{definition} {\rm {\cite{Ger} \cite[Definition 7.5.1]{Mac}}} Let
  $R$ be a ring with unit and let $C $ be an $R$-algebra. Let
  $\mathfrak{a}=\oplus_{n\in \mathbb{Z}}\mathfrak{a}^n $ be a graded
  algebra over $C$.  $\mathfrak{a}$ is a \emph{Gerstenhaber algebra}
  if there is an associative product $\wedge$ and a graded commutative
  product $[ {-} \bullet {-} ]$ satisfying the following axioms. When
  $a\in \mathfrak{a}^n$, let $|a|$ denote its degree $n$. For $a\in
  \mathfrak{a}^{|a|}$, $b\in \mathfrak{a}^{|b|}$, $c\in
  \mathfrak{a}^{|c|}$,
\begin{eqnarray}
 a\wedge b\in \mathfrak{a}^{|a|+|b|}, && b\wedge
a=(-1)^{|a||b|}a\wedge b. \\
 {[ {a} \bullet {b} ]} \in \mathfrak{a}^{|a|+|b|-1},
 && [ {a} \bullet {b} ]=-(-1)^{(|a|+1)(|b|+1)}[ {b} \bullet {a} ].
\label{commutative}
\end{eqnarray}
 \begin{equation}\label{jacobi} (-1)^{(|a|+1)(|c|+1)}[ {[ {a}
\bullet {b} ]} \bullet {c} ] +(-1)^{(|b|+1)(|a|+1)}[ {[ {b} \bullet
{c} ]} \bullet {a} ] +(-1)^{(|c|+1)(|b|+1)}[ {[ {c} \bullet {a} ]}
\bullet {b} ]=0.
\end{equation}
\begin{equation}  [ {a} \bullet {b\wedge c} ]=[ {a} \bullet {b} ]\wedge
c+(-1)^{(|a|+1)|b|}b\wedge [ {a} \bullet {c} ]. \label{distributive}
\end{equation}
\end{definition}
On the other hand, we have the following construction.
\begin{definition} A differential graded algebra  is a graded algebra $%
\mathfrak{a}=\oplus_{n\in \mathbb{Z}}\mathfrak{a}^n$ with a graded
commutative product $\wedge$ and a differential $d$ of degree $+1$,
i.e. a map $d:\mathfrak{a}\to \mathfrak{a}$ such that
\begin{equation}
 d(\mathfrak{a}^n)\subseteq \mathfrak{a}^{n+1}, \quad d\circ d=0,
 \quad
 d(a\wedge b)=da\wedge b+(-1)^{|a|}a\wedge db. \label{compatible
 with wedge}
\end{equation}
\end{definition}

\begin{definition}
  Let $\mathfrak{a}=\oplus_{n\in \mathbb{Z}}\mathfrak{a}^n$ be a
  graded algebra over $C$ such that $(\mathfrak{a}, [ {-} \bullet {-}
  ], \wedge)$ form a Gerstenhaber algebra and $(\mathfrak{a}, \wedge,
  d)$ form a differential graded algebra. If in addition
\begin{equation}
d[ {a} \bullet {b} ]=[ {da} \bullet {b} ]+(-1)^{|a|+1}[ {a} \bullet
{db} ], \label{compatible with br}
\end{equation}
for all $a$ and $b$ in $\mathfrak{a}$, then $(\mathfrak{a}, [ {-}
\bullet {-} ], \wedge, d)$ is a differential Gerstenhaber algebra
{\rm (DGA)}.
\end{definition}

For any Gerstenhaber algebra, $\mathfrak{a}^1$ with the induced
bracket is a Lie algebra.  Conversely, suppose $\mathfrak{a}^1$ is a
finite dimensional algebra over the complex or real numbers,
equipped with a differential compatible with the Lie bracket.  Then a
straightforward induction allows one to construct a DGA structure on
the exterior algebra of $\mathfrak{a}^1$.

\begin{lemma}\label{natural extension} Let $\mathfrak{a}^1$ be a finite
  dimensional Lie algebra with bracket $[ {-} \bullet {-} ]$. Let
  $\mathfrak{a} $ be the exterior algebra generated by
  $\mathfrak{a}^1$. Then the Lie bracket on $\mathfrak{a}^1$ uniquely
  extends to a bracket on $\mathfrak{a}$ so that $(\mathfrak{a}, [ {-}
  \bullet {-} ], \wedge)$ is a Gerstenhaber algebra.

  If, furthermore, an operator $d:\mathfrak{a}^1\to\mathfrak{a}^2$ is
  extended as in {\rm (\ref{compatible with wedge})}, then
  $(\mathfrak{a}, [ {-} \bullet {-} ], \wedge)$ is a differential
  Gerstenhaber algebra if and only if
\begin{equation}
d[ {a} \bullet {b} ]=[ {da} \bullet {b} ]+[ {a} \bullet {db} ],
\label{compatible with br one degree 1}
\end{equation}
for all $a$ and $b$ in $\mathfrak{a}^1$.
\end{lemma}

\begin{definition}
  A homomorphism of differential graded Lie algebras is called a
  quasi-isomorphism if the map induced on the associated cohomology
  groups is a linear isomorphism.

  A quasi-isomorphism of differential Gerstenhaber algebras is a
  homomorphism of DGAs that descends to an isomorphism of cohomology
  groups.
\end{definition}
Note that in the latter case the isomorphism is one of Gerstenhaber
algebras.

\subsection{DGA of complex structures}

Suppose $J$ is an \emph{integrable complex structure} on
$\mathfrak{g}$.  i.e. $J$ is an endomorphism of $\mathfrak{g}$ such
that $J\circ J=-1$ and
\begin{equation}  \label{integrable}
  [ {x} \bullet {y} ]+J[ {Jx}
\bullet {y} ]+J[ {x} \bullet {Jy} ]-[ {Jx} \bullet {Jy} ]=0.
\end{equation}
Then the $\pm i$ eigenspaces $\mathfrak{g}^{(1,0)}$ and
$\mathfrak{g}^{(0,1)} $ are complex Lie subalgebras of the
complexified algebra $\mathfrak{g}_\CC$. Let $\mathfrak{f}$ be the
exterior algebra generated by $\mathfrak{g}
^{(1,0)}\oplus\mathfrak{g}^{*(0,1)}$, i.e.
\begin{equation}
\mathfrak{f}^n:=\wedge^n(\mathfrak{g}^{(1,0)}\oplus\mathfrak{g}^{*(0,1)}),
\quad \mbox{ and } \quad \mathfrak{f}=\oplus_n\mathfrak{f}^n.
\end{equation}
The integrability condition in (\ref{integrable}) implies that
$\mathfrak{f}^1$ is closed under the \emph{Courant bracket}
\begin{equation}
[ x+\alpha \bullet y+\beta ]:=[x, y]+\iota_xd\beta-\iota_yd\alpha.
\end{equation}
A similar construction holds for the conjugate
${\overline{\mathfrak{f}}}$, generated by $
\mathfrak{g}^{(0,1)}\oplus\mathfrak{g}^{*(1,0)}.  $

Recall that if $(\mathfrak{g}, [ {-} \bullet {-} ])$ is a Lie algebra,
the Chevalley-Eilenberg (C-E) differential $d$ is defined on the dual
vector space $\mathfrak{g}^*$ by the relation
\begin{equation}  \label{C-E}
d\alpha (x,y):=-\alpha ([ {x} \bullet {y} ]),
\end{equation}
for $\alpha \in \mathfrak{g}^*$ and $x,y\in \mathfrak{g}$.  This
operator is extended to the exterior algebra $\wedge\mathfrak{g}^*$ by
derivation.  The identity $d\circ d=0$ is equivalent to the Jacobi
identity for the Lie bracket $[ {-} \bullet {-} ]$ on $\mathfrak{g}$.
It follows that $(\wedge \mathfrak{g}^*, d)$ is a differential graded
algebra.

The natural pairing on $(\lie g\oplus\lie g^*)\otimes \mathbb C$,
induces a complex linear isomorphism $ (\mathfrak{f}^1)^*\cong
{\overline{\mathfrak{f}}}^1$.  Therefore, the C-E differential of the
Lie algebra ${\overline{\mathfrak{f}}}^1$ is a map from
$\mathfrak{f}^1$ to $\mathfrak{f}^2$. Denote this operator by
$\overline\partial$.  Similarly, we denote the C-E differential of $
\mathfrak{f}^1$ by $\partial$. It is well known that the maps
\begin{equation}
\overline\partial: \mathfrak{g}^{*(0,1)}\to
\wedge^2\mathfrak{g}^{*(0,1)}, \quad \mbox{ and } \quad
\overline\partial: \mathfrak{g}^{(1,0)}\to
\mathfrak{g}^{(1,0)}\otimes \mathfrak{g}^{*(0,1)}
\end{equation}
are respectively given by
\begin{equation}
\overline\partial\overline\omega=(d\overline\omega)^{0,2}, \quad
(\overline\partial T){\overline W}= [ {{\overline W}} \bullet
{T}]^{1,0}
\end{equation}
for any $\omega$ in $\mathfrak{g}^{*(0,1)}$, $T\in
\mathfrak{g}^{1,0}$ and ${\overline W}\in \mathfrak{g}^{0,1}$.

If $\{T_\ell: 1\leq\ell\leq n\}$ forms a basis for $
\mathfrak{g}^{1,0}$ and $\{\omega^\ell:1\leq \ell\leq n\}$ the
dual basis in $\mathfrak{g}^{*(1,0)}$, then we have
\begin{equation}
\overline\partial{\overline\omega}^\ell=(d{\overline\omega}^\ell)^{0,2},
\quad (\overline\partial T)=\sum_\ell\bom{\ell}\wedge [ {{\overline
T}_\ell} \bullet T_j]^{1,0}
\end{equation}

Based on Lemma \ref{natural extension}, it is an elementary
computation to verify that the  quadruples $(\mathfrak{f}, [ {-}
\bullet {-} ], \wedge, \overline\partial)$ and
$(\overline{\mathfrak{f}}, [ {-} \bullet {-} ], \wedge, \partial)$
are differential Gerstenhaber algebras.

For a given Lie algebra $\lie g$ and a choice of invariant complex
structure $J$, we denote the differential Gerstenhaber algebra
$(\mathfrak{f}, [ {-} \bullet {-} ], \wedge, \overline\partial)$ by
$\DGA ({\lie g}, J)$.

The following observation relying on the nature of the
$\overline\partial$ and $\partial$ will be helpful, although
apparently obvious.

\begin{lemma}\label{duality} Given a complex linear identification
  ${\overline{\lie f}}^1\cong ({\lie f}^1)^*$, the Lie algebras
  $({\lie f}^1, [-\bullet -])$, $({\overline{\lie f}}^1, [-\bullet
  -])$ and the graded differential algebras $({\lie f}, \dbar)$,
  $({\overline{\lie f}},
\partial)$ determine each other.
\end{lemma}

\subsection{DGA of symplectic structures}\label{symp}
Let $\mathfrak{h}$ be a Lie algebra over $\RR$. The exterior algebra
of the dual $\mathfrak{h}^*$ with the C-E differential $d$ is a
differential graded Lie algebra.

Suppose that $O:\mathfrak{h}\to\mathfrak{h}^*$ is a real linear map.
Define a bracket $[ {-} \bullet {-} ]_O$ on $\mathfrak{h}^*$ by
\begin{equation}
[ {\alpha} \bullet {\beta} ]_O:=O[ {O^{-1}\alpha} \bullet
{O^{-1}\beta} ].
\end{equation}
It is a tautology that $(\mathfrak{h}^*, [ {-} \bullet {-} ]_O)$
becomes a Lie algebra, with the map $O$ understood as a Lie algebra
homomorphism.

\begin{definition}
A linear map $O:\mathfrak{h}\to\mathfrak{h}^*$ from a Lie algebra to
its dual is said to be compatible with the C-E differential if  for
any $\alpha, \beta$ in $\mathfrak{h}^*$,
\begin{equation}
d[ {\alpha} \bullet {\beta} ]_O=[ {d\alpha} \bullet {\beta} ]_O+[
{\alpha} \bullet {d\beta} ]_O.
\end{equation}
\end{definition}

Due to Lemma \ref{natural extension}, the next observation is a
matter of definitions.

\begin{lemma}\label{construction}
  Suppose $\mathfrak{h}$ is a Lie algebra, and take an element $O$ in
  $\Hom(\mathfrak{h}, \mathfrak{h}^*)$ compatible with the C-E
  differential. Then $(\wedge^\bullet\mathfrak{h}^*, [ {-} \bullet {-}
  ]_O, \wedge, d)$ is a differential Gerstenhaber algebra.
\end{lemma}

When the algebra $\lie h$ has a symplectic form $\Omega$, the
contraction with $\Omega$ defines an $O$ as in the above lemma. In
such case, the differential Gerstenhaber algebra $(\wedge^\bullet{\lie
  h}^*, [-\bullet -]_\Omega, \wedge, d)$ {\it after complexification}
is denoted by $\DGA ({\lie h}, \Omega)$.

\section{Complex Structures on Nilpotent Algebras}\label{sec:nil}

\subsection{General Theory}
Let $\lie g$ be a Lie algebra over $\mathbb R$ or $\mathbb C$. The
\emph{lower central series} of $\lie g$ is the sequence of
subalgebras $\lie g_{p+1} \subset \lie g_p\subset\lie g$ given by
\[
\lie g_0=\lie g, \quad \lie g_p = [\lie g_{p-1}\bullet\lie g]. \] A Lie
algebra $\lie g$ is $s$-step nilpotent if $s$ is the smallest
integer such that $\lie g_s=\{0\}$.  Defining $V_p$ to be the
annihilator of $\lie g_p$ one has \emph{the dual sequence} $\lie
g^*\supset V_{p}\supset V_{p+1}$.  The dual sequence may also be
defined recursively as
\[ V_0 = \{0\}, \quad V_p = \{
\alpha\in\lie{g}^*:d\alpha\in\Lambda^2 V_{p-1}\}. \]
We note that if the subscript $\CC$ denotes complexification of vector
spaces and Lie algebras, then $V_p(\lie g)_\CC = V_p(\lie g_\CC)$.
Write $n_p=\dim V_p$.  A \emph{Malcev} basis for $\lie g^*$ is a basis
chosen such that $e^1,\dots,e^{n_1}$ is a basis for $V_1$,
supplemented with $e^{n_1+1},\dots,e^{n_2}$ to form a basis for $V_2$,
et cetera.  For such a basis one has $de^p\in\Lambda^2\langle
e^1,\dots, e^{p-1}\rangle$. The short-hand notation
$12:=e^{12}:=e^1\wedge e^2$ is convenient.  Using this one may
identify a Lie algebra by listing its structure equations with respect
to a Malcev basis as $(de^1,\dots,de^n)$.  For instance we may write
$\lie g=(0,0,a{12})$ to mean the Lie algebra $\lie g$ generated by the
relations $de^1=0=de^2$, $de^3=a e^1\wedge e^2$.  This has the single
non-trivial bracket $[e_1\bullet e_2]=-ae_3$.

\begin{lemma}\label{quasi}
  Suppose $\lie h$ and $\lie k$ are Lie algebras, $\lie h$ is
  nilpotent and $\phi\colon(\bw \lie h^*,d)\to(\bw\lie k^*,d)$ is a
  quasi-isomorphism of the associated differential graded algebras.
  Then $\phi$ is an isomorphism.
\end{lemma}
  \bproof
  If $\phi\colon\lie g^*\to\lie h^*$ is a homomorphism of the
  associated differential graded algebras then $\phi(V_p(\lie
  g))\subset V_p(\lie h)$.  When $\phi$ is furthermore a
  quasi-isomorphism then the restriction $V_1(\lie g)\to V_1(\lie h)$
  is an isomorphism of vector spaces.  Suppose that $\phi$ restricted
  to $V_{p-1}(\lie g)$ is an isomorphism onto $V_{p-1}(\lie h)$.  Then
  clearly the induced map $\Lambda^2V_{p-1}(\lie
  g)\to\Lambda^2V_{p-1}(\lie h)$ is also an isomorphism.  Suppose that
  $a\in V_p(\lie g)$ satisfies $\phi(a)=0$.  Then
  $da\in\Lambda^2V_{p-1}(\lie g)$ satisfies $\phi(da)=0$.  But then
  $da=0$, so $a\in V_1$, and $\phi(a)=0$ actually implies $a=0$.
  \eproof

\begin{proposition}\label{key technical}
  Suppose that $\lie g$ and $\lie h$ are finite dimensional nilpotent
  Lie algebras, $J$ is an integrable complex structure on $\lie g$ and
  $O\colon\lie h\to\lie h^*$ is a linear map compatible with the C-E
  differential on $\lie h$.  Then a homomorphism $\phi$ from $\DGA
  ({\lie g}, J)$ to $\DGA ({\lie h}, O)$ is a quasi-isomorphism if and
  only if it is an isomorphism.
\end{proposition}
\bproof As a quasi-isomorphism of DGAs, $\phi$ is a
quasi-isomorphism of the underlying exterior differential algebras:
\[ \phi: (\wedge^*{\lie f}^1, \wedge, \dbar)\to (\wedge^*{\lie h}_\CC^*,
\wedge, d).
\]
The last lemma shows that it has to be an isomorphism. \eproof

Two special types of complex DGAs were introduced above: Those coming
from an integrable complex structure $J$ on a real algebra $\lie g$,
denoted $\DGA(\lie g,J)$ and those derived from a linear
identification $O\colon\lie h\to\lie h^*$ compatible with the
differential.  For the latter we write $\DGA(\lie h,O)$.  A problem
related to ``weak mirror symmetry'' is: given an algebra $\lie g$ and
a complex structure $J$, when does an $\lie h$ and $O$ exist so that
$\DGA(\lie g,J)$ is quasi-isomorphic to $\DGA(\lie h,O)$?  For
nilpotent algebras we shall see below that this always is the case.

Given $J$ on $\lie g$ write $\lie k$ for the Lie algebra ${\bar{\lie
    f}}^1$ consisting of degree one elements in $\DGA(\lie g,J)$.  As
$\lie h$ is complex we may speak of the complex conjugate algebra
consisting of the conjugate vector space $\bar{\lie h}$.  Let
$c\colon\lie h\to\bar{\lie h}$ be the canonical map such that
$c(ax)=\bar a x$ for complex $a$ and $x$ in $\lie h$.  Then
$[x,y]_c:=c[c(x),c(y)]$ equips $\bar{\lie h}$ with a Lie bracket.  We
say that $\lie h$ is \emph{self-conjugate} if a complex linear
isomorphism $\lie h\to \bar{\lie h}$ exists.

\begin{proposition}\label{second technical}
  Let $\lie g$ be a Lie algebra with complex structure $J$. Let
  $\DGA(\lie g, J)$ be its differential Gerstenhaber algebra. Write
  $\lie h$ for the Lie algebra $\lie f^1$ and suppose that $\lie h$ is
  self-conjugate.  Then there exists a complex linear isomorphism
  $O\colon\lie h\to{\lie h}^*$ compatible with the C-E
  differential $d$ on $\lie h$ so that $\DGA(\lie h,O)$ is isomorphic
  to $\DGA(\lie g, J)$.
\end{proposition}

\bproof We construct the map $O$.  Let $\phi \colon \lie h \to \lie
f^1$ be the identification of $\lie h$ as the Lie algebra given by
$\lie f^1$.  Composing on both sides with complex conjugation gives
the isomorphism $\bar\phi := c\circ\phi\circ c \colon {\bar{\lie h}}
\to {\bar{\lie f}^1}$ of Lie algebras.  Taking the identifications
${\bar{\lie f}^1}\cong(\lie f^1)^*$ and $\bar{\lie h}\cong\lie h$ into
account gives the isomorphism
\begin{equation}
\psi:  \lie h\cong{\bar{\lie
    h}} \stackrel{{\bar\phi}}{\rightarrow} {\bar{\lie
    f}^1}\cong(\lie f^1)^*
\end{equation}
of Lie algebras.  The dual map $\psi^*$ is now an isomorphism of
exterior differential algebras
\begin{equation}\label{natural1}
\psi^*: \wedge^*{\lie f}^1 \to \wedge^*{\lie h}^*, \quad \psi^*\circ \dbar
=d \circ \psi^*.
\end{equation}
We claim that the following composition
\begin{equation}
O: \lie h \stackrel{\phi}{\rightarrow} {\lie f}^1
 \stackrel{{\psi}^*}{\rightarrow} \lie h^*
\end{equation}
is compatible with $d$. By Lemma~\ref{natural extension} this is the
case if equation \eqref{compatible with br one degree 1} holds for the
bracket $\sbr{\alpha}{\beta}_O=:O\sbr{O^{-1}\alpha}{O^{-1}\beta}$.  But
\begin{equation}\label{natural2}
\sbr{\alpha}{\beta}_O= \psi^*\circ \phi
\sbr{\phi^{-1}\circ
(\psi^*)^{-1}\alpha}{\phi^{-1}\circ
(\psi^*)^{-1} \beta}= \psi^*
\sbr{ (\psi^*)^{-1}\alpha}{ (\psi^*)^{-1} \beta}.
\end{equation}
Hence,
\begin{eqnarray*}
 d\sbr{\alpha}{\beta}_O&=&d(\psi^*
\sbr{ (\psi^*)^{-1}\alpha}{ (\psi^*)^{-1} \beta} ) \\
&=&\psi^* (\dbar \sbr{ (\psi^*)^{-1}\alpha}{ (\psi^*)^{-1} \beta} )\\
&=&\psi^* (\sbr{ \dbar (\psi^*)^{-1}\alpha}{ (\psi^*)^{-1} \beta} )+\psi^* ( \sbr{ (\psi^*)^{-1}\alpha}{ \dbar (\psi^*)^{-1} \beta} ) )\\
&=&\psi^* (\sbr{ (\psi^*)^{-1}d\alpha}{ (\psi^*)^{-1} \beta} )+\psi^* ( \sbr{ (\psi^*)^{-1}\alpha}{ (\psi^*)^{-1} d\beta} ) )\\
&=& \sbr{ d\alpha}{ \beta}_O+\sbr{\alpha}{d \beta}_O.
\end{eqnarray*}
By Lemma \ref{construction}, $\DGA(\lie k, O):=(\wedge^*{\lie k}^*,
\sbr{-}{-}_O, \wedge, d)$ forms a differential Gerstenhaber algebra.
It is clear from (\ref{natural1}) and (\ref{natural2}) that the map
$\psi^*$ yields an isomorphism from $\DGA(\lie g, J)$ to $\DGA(\lie k,
O)$.  \eproof

It should be noted that the map $O$ is not necessarily skew-symmetric,
nor is it automatically closed when it is skew.  In particular the DGA
structure obtained above does not necessarily arise from contraction
with a symplectic structure.  Also note that the condition ${\bar{\lie
    k}}\cong\lie k$ is satisfied precisely when $\lie k$ is the
complexification of some real algebra.  Whilst in the context of
six-dimensional nilpotent algebras this is always the case, there
exist non-isomorphic real algebras having the same complexification.

\subsection{Nilpotent complex structures}
 An almost complex structure $J$ on $\lie g$ may be
given by a choice of basis $\omega=\{ \omega^{k}, 1\leq k \leq m
\}$, $2m=\dim_\RR\lie g$, of the space of $(1,0)$-forms in the
complexified dual $\lie g_\CC^*$. Such a basis may equivalently be
given as a basis $e=(e^1,\dots,e^{2m})$ of $\lie g^*$ so that
$e^2=Je^1$, or $\omega^1=e^1+ie^2$, and so on. When $e$ and $\omega$
are related in this way we will write $e=e(\omega)$ or
$\omega=\omega(e)$.  The almost complex structure is then
\emph{integrable} or simply a \emph{complex structure} if the ideal
in $\Lambda^*\lie g^*_\CC$ generated by the $(1,0)$-forms is closed
under exterior differentiation.  For a nilpotent Lie algebra, an
almost complex structure is integrable if there exists a basis
$(\om{j})$ of $(1,0)$-forms so that $d\om1=0$ and for $j>1$,
$d\om{j}$ lies in the ideal generated by $\om1,\dots,\om{j-1}$.
Equivalently,
\begin{equation}
  0 = d(\om1\wedge\om2\wedge \dots\wedge
  \om{p}),\quad p=1,\dots,m.
\end{equation}
Let the set of such bases be denoted $\Omega(\lie g,J)$.

On nilpotent Lie algebras certain complex structures are
distinguished.  Among these are complex structures such that
$[X,JY]=J[X,Y]$.  Equivalently,
$d\om{p}\in\Lambda^2\langle\om1,\dots,\om{p-1}\rangle$.  These are the
complex structures for which $\lie g$ is the real algebra underlying a
complex Lie algebra.  At the opposite end to these are the \emph{abelian
  complex structures} which satisfy $[JX,JY]=[X,Y]$ \cite{BDM}.
Equivalently, the $+i$-eigenspace of $J$ in $\lie g_\CC$ is an abelian
subalgebra of $\lie g_\CC$.  In particular abelian $J$s are always
integrable.  In terms of $(1,0)$-forms a complex structure is abelian
if and only if there exists an $\omega$ in $\Omega(\lie g,J)$ such that
$d\om{j}$ is in the intersection of the two ideals generated by
$\om1,\dots,\om{j-1}$ and $\bom1,\dots,\bom{j-1}$, respectively.

The concept of abelian complex structures may be generalized to that
of \it nilpotent \rm complex structures \cite{CFGU}.  A nilpotent
almost complex structure may be defined as an almost complex
structure with a basis of $(1,0)$-forms such that
\begin{equation}
  d\om{p}\in\Lambda^2\langle\om1,\dots,\om{p-1},\bom1,\dots,\bom{p-1}\rangle.
\end{equation}
For a given algebra $\lie g$ and nilpotent almost complex structure
$J$ we write $P(\lie g,J)$ for the set of such bases.  Nilpotent
almost complex structures are not necessarily integrable.  If a
nilpotent $J$ is integrable, then  $P(\lie g,J)\subset\Omega(\lie
g,J)$.  A nilpotent complex structure  is abelian if and only if
\begin{equation}
  0 = d(\bom1\wedge\bom2\wedge \dots\wedge\bom{p-1}\wedge
  \om{p}), \quad p=1,\dots,m.
\end{equation}
It is apparent that abelian complex structures are nilpotent.

In subsequent presentation, we suppress the wedge product sign.

\subsection{Six-dimensional algebras}
Some of the results of this section may be regarded as a re-organization of
past results in terms of invariants relevant to our further
analysis.  Our key references are \cite{Salamon} and \cite{Ugarte}.
To name specific isomorphism classes of six-dimensional nilpotent
Lie algebras, we use the notation $\lie h_n$ as given in
\cite{CFGU}.

Suppose then $\dim_\RR\lie g=6$. Let $J$ be a nilpotent almost complex
structure on $\lie g$.  The \emph{structure equations} for an
\emph{integrable} element $\omega$ in $P(\lie g,J)$ are \cite{CFGU}
\begin{gather}\label{eq:6'}
  \begin{cases}
    d\om1 = 0,\\
    d\om2 = \epsilon\om1\bom{1},\\
    d\om3 = \rho\om1\om2 + A\om1\bom1 + B\om1\bom2 +
    C\om2\bom1 + D\om{2}\bom2.
  \end{cases}
\end{gather}
for complex numbers $\epsilon,\rho,A,B,C,D$.  Note that $dd\om3=0$
forces $D\epsilon=0$.  Moreover, if $\epsilon$ is not zero, $\om3$
may be replaced with $\epsilon\om3-A\om2$ so after re-scaling the
$\om{j}$ one obtains the \emph{reduced structure
  equations} \cite{Ugarte}
\begin{gather}\label{eq:6}
  \begin{cases}
    d\om1 = 0,\\
    d\om2 = \epsilon\om1\bom{1},\\
    d\om3 = \rho\om1\om2 + (1-\epsilon)A\om1\bom1 + B\om1\bom2 +
    C\om2\bom1 + (1-\epsilon)D\om{2}\bom2,
  \end{cases}
\end{gather}
where $\epsilon$ and $\rho$ are either $0$ or $1$ and $A,B,C,D$ are
complex numbers.

To avoid ambiguity we rule out the case $\epsilon\not=0$, $d\om3=0$
for any form of the structure equations as this is equivalent to
$\epsilon=0, ~d\om3=\om{1}\bom1$.

Given structure equations~\eqref{eq:6'} for a nilpotent complex
structure, we will calculate $\DGA(\lie g,J)$  in
Section~\ref{sec:isom-class-lie}.  However, if we take~\eqref{eq:6'}
as a starting point, it is not obvious to recognize the real algebra
$\lie g$ which underlies the complex structure. We shall first
provide a way to do this that will fit the purpose of this paper.

For this task, we identify invariants of $P({\lie g}, J)$.
The most immediate invariants are the dimensions of the vector spaces
in the dual sequence $V_0\subset V_1\subset\dots$ for $\lie g_\CC$. As
the inclusions
\begin{equation}
  \label{eq:34}
  V_1\supset\langle\om1,\bom1,\om2+\bom2\rangle, \qquad
  V_2\supset\langle\om1,\bom1,\om2,\bom2\rangle
\end{equation}
always hold, $V_3=\lie g_\CC^*$ for any $6$-dimensional nilpotent
algebra with nilpotent complex structure.  Define
\begin{equation}
  n=(n_1,n_2)=(\dim V_1, \dim V_2). \end{equation} We now  collect
several facts on these particular invariants.

\begin{lemma}\label{dependence} Given a nilpotent complex structure $J$ on a six-dimensional nilpotent algebra
$\lie g$, the following hold:
\begin{enumerate}[{\rm(a)}]
\item $3\leq n_1\leq 6$, $4\leq n_2\leq 6$ and $n_1\leq n_2$.
\item There exists $\omega$ in $P({\lie g}, J)$ such that $\epsilon=0$
or $\epsilon=1$.
\item If $\epsilon=1$, there exists $\omega$ in $P({\lie g}, J)$ such
that $A=D=0$.
\item If $\epsilon=0$, then $n_2=6$.\label{item:1}
\item $\rho=0$ if and only if $J$ is an abelian complex structure.
\item Let $d$ be the dimension of the complex linear span of $d\om3$
  and $d\bom3$. Then $d\leq 1$ if and only if
  \begin{equation}
    \label{eq:32}
  \rho=0,\qquad  |B|^2=|C|^2,\qquad A\bar D=\bar A D,\qquad A\bar B=\bar A
    C,\qquad D\bar B = \bar D C.
  \end{equation}\label{item:14}
  \end{enumerate}
\end{lemma}

Based on the above information, we  re-organize some of the data
from \cite[Theorem 2.9]{Ugarte} and \cite[Table A.1]{Salamon}.

\begin{lemma}\label{lem:3}
  Suppose a complex structure on $\lie
  g$ is given with structure constants as in~\eqref{eq:6'} with
  $\epsilon\in\{0,1\}$.
  \begin{enumerate}[{\rm(a)}]
  \item $n=(6,6)$ if and only if $\lie g\cong\lie h_1=(0,0,0,0,0,0)$.
  \item $n=(5,6)$ if and only if $\epsilon=0$ and $d=1$.  In this case,
    \begin{gather*}
      \lie g\cong
      \begin{cases}
        {\lie h_{8}}=(0,0,0,0,0,12),\\
        \lie h_{3}=(0,0,0,0,0,12+34).
      \end{cases}
    \end{gather*}
  \item If $n=(4,6)$, then $\epsilon=0$ and $d=2$. The Lie algebra is
    \begin{gather*}
      \lie g\cong
      \begin{cases}
        \lie h_{6}=(0,0,0,0,12,13),\\
        \lie h_{2}=(0,0,0,0,12,34),\\
        \lie h_{4}=(0,0,0,0,12,13+42),\\
        \lie h_{5}=(0,0,0,0,13+42,14+23).
      \end{cases}
    \end{gather*} \label{item:d}
  \item If $n=(3,6)$, then $\epsilon=1, \rho\neq 0$, and there exists an element $\sigma$ in $P({\lie g}, J)$
  such that  $d\te3 =
    \te1(\te2+\bte2)$. Moreover, $\lie g\cong\lie
    h_7=(0,0,0,12,13,23)$.
    \label{item:2}
  \item If $n=(4,5)$, there exists $\sigma$ in $P({\lie g}, J)$ such that $d\te3 =
    \te1\bte2 + \te2\bte1$.  The structure equations for $e(\sigma)$
    are $(0,0,0,12,0,14-23)$ and so $\lie g\cong\lie h_{9}=(0,0,0,0,12,14+25)$. \label{item:3}
  \item \label{item:4} If $n=(3,5)$, then there exists $\sigma$ in $P({\lie
      g}, J)$ such that $d\te3 = (B-\bar C)\te1\te2+B\te1\bte2 +
    C\te2\bte1$.  For all such $\sigma$, $(B-\bar C)$ is non-zero.
    Moreover,
    \begin{gather*}
      \lie g\cong
      \begin{cases}
        {\lie h_{10}}=(0,0,0,12,13,14),\\
        {\lie h_{12}}=(0,0,0,12,13,24),\\
        {\lie h_{11}}=(0,0,0,12,13,14+23).
      \end{cases}
    \end{gather*}
  \item If $n=(3,4)$, then
    \begin{gather*}
      \lie g\cong
      \begin{cases}
        \lie h_{16}=(0,0,0,12,14,24),\\
        \lie h_{13}=(0,0,0,12,13+14,24),\\
        \lie h_{14}=(0,0,0,12,14,13+42),\\
        \lie h_{15}=(0,0,0,12,13+42,14+23).
      \end{cases}
    \end{gather*}
  \end{enumerate}
\end{lemma}

\bproof
  The statements for $n=(6,6)$ and $(5,6)$
  are elementary.

  When $n=(4,6)$ the cases listed in (\ref{item:d}) are the only
  possibilities given by the classification of \cite{Ugarte}.

  When $n=(3,6)$, then $V_1=\langle\om1,\bom1,\om2+\bom2\rangle$ and
  $d\om3\in\Lambda^2V_1$.  It follows that $d\om3=\rho\om1(\om2+\bom2)$
  and so we have~\eqref{item:2}.

  If $n_2=5$, then $\epsilon=1$, and we may take $A=D=0$ as noted in
  the previous lemma.  By~\eqref{eq:34}, a complex number $u\not=0$
  exists so that $ud\om3+\bar ud\bom3\in\Lambda^2 V_1$.

  If in addition $n_1=3$, then $V_1=\langle \om1, \bom1,
  \om2+\bom2\rangle$.  Taking
  $d\om3=\rho\om1\om2+B\om1\bom2+C\om2\bom1$ gives $u\rho=uB-\bar
  u\bar C$.  Now setting $\te1=\om1,~\te2=\om2$ and $\te3=u\om3$ puts
  the structure equations in the form (\ref{item:4}).  Note that if
  $B=\bar C$ in (\ref{item:4}) then $d\te3$ and $d\bte3$ are linearly
  dependent and so $n_1=4$.

  If, on the other hand, $n_1=4$ then $V_1=\langle \om1, \bom1,
  \om2+\bom2, \om3+\lambda\bom3\rangle$ for some $\lambda$.  Then
  $\rho=0$, $B=\lambda {\overline C}$ and $C=\lambda {\overline B}$,
  since $d\om3+\lambda d\bom3=0$. In particular,
  $d\om3=B\om1\bom2+\lambda {\overline B}\om2\bom1$. This yields case
  (\ref{item:3}).

  The remaining case is $n=(3,4)$. Since this is the minimum
  possible combination for the invariant $n$, by exclusion all remaining
  nilpotent complex structures found in \cite{Salamon} and
  \cite{Ugarte} are covered in this case.  \eproof

\begin{corollary}\label{cor:1}
  Let $J$ be a nilpotent complex structure on a nilpotent Lie algebra
  $\lie g$.  Then complex structure is abelian if
   $n=(6,6), (5,6), (4,5)$.    It is not abelian if $n_1=3$ and $n_2>4$.
\end{corollary}

\subsection{More invariants of nilpotent complex structures}
Given $\omega$ in  $P(\lie g,J)$. Suppose that its structure
equations are (\ref{eq:6'}). Let $\sigma$ be another element in
$P(\lie g,J)$. Viewing $\sigma$ and $\omega$ as row vectors, then
$\sigma=(\te1,\te2, \te3)$ and $\omega=(\om1,\om2,\om3)$ are related
by a matrix:  $\te{j}=\sigma^j_k\om{k}$. This must be of the form
\begin{gather}\label{eq:29}
  \sigma(\omega): =
  \begin{pmatrix}
    \sigma^1_1&\sigma^2_1&\sigma^3_1\\
    \sigma^1_2&\sigma^2_2&\sigma^3_2\\
    0&0&\sigma^3_3
  \end{pmatrix},
\end{gather}
with $\epsilon\sigma^1_2=0$.  So when $\epsilon\not=0$ the matrix
$\sigma(\omega)$ is upper triangular.  Write $\Delta(\sigma,\omega)$
for the determinant of the transformation $\sigma(\omega)$, so that
$\te1\te2\te3 = \Delta(\sigma,\omega)\om1\om2\om3$, and
$\Delta(\sigma,\omega)^{-1} = \Delta(\omega,\sigma)$. Define
$\Delta'(\sigma,\omega)$ by $\te1\te2 =
\Delta'(\sigma,\omega)\om1\om2$ so that $\Delta(\sigma,\omega) =
\sigma^3_3\Delta'(\sigma,\omega)$.  The space $\cal A$ of matrices as
in~\eqref{eq:29} may be considered the automorphism group of the
nilpotent complex structure, and $P(\lie g,J)$ is the orbit of
$\omega$ under the multiplication of elements in $\cal A$.

Consider the two functions $\Delta_1\colon P(\lie g,J)\to\mathbb C$,
$\Delta_2\colon P(\lie g,J)\to\mathbb R$ defined respectively by
\begin{gather}
  d\te3\wedge d\te3 = 2\Delta_1(\sigma)\te1\bte1\te2\bte2,\\
  d\te3\wedge d\bte3 = 2\Delta_2(\sigma)\te1\bte1\te2\bte2.
\end{gather}
In terms of the structure constants for $\omega$,
\begin{align}
  \Delta_1(\omega) &= AD - BC,\label{eq:35}\\
  \Delta_2(\omega) &= \tfrac12\left[{|B|}^2+{|C|}^2-A\bar D-\bar A
    D
    -{|\rho |}^2\right]\label{eq:36}.
\end{align}
If $\sigma=\sigma(\omega)$ then
\begin{align*}
  d\te3\wedge d\te3 &= (\sigma^3_3)d\om3\wedge d\om3 = (\sigma^3_3)\Delta_1(\omega)\om1\bom1\om2\bom2 \\
  &= (\sigma^3_3)^2
  \Delta_1(\omega){|\Delta'(\omega,\sigma)|}^2 \te1\bte1\te2\bte2.
\end{align*}
Therefore
\begin{gather}
  \Delta_1(\sigma) =
  \Delta_1(\omega){|\Delta'(\omega,\sigma)|}^2(\sigma^3_3)^2,
  \intertext{and similarly}
  \Delta_2(\sigma) =
  \Delta_2(\omega) {|\Delta'(\omega,\sigma)|}^2 {|\sigma^3_3|}^2.
\end{gather}
By choosing $\sigma$ appropriately we may assume that $\Delta_1$ is
either $0$ or $1$.  We observe that if $\Delta_1$ is non-zero in some
basis then it is non-zero in every basis. In this situation
$\Delta_2/{|\Delta_1|}$ is invariant under transformations of the
form~\eqref{eq:29}.  Note that $\Delta^2_2-{|\Delta_1|}^2$ is scaled
by a positive constant by an automorphism, so the sign of $\Delta^2_2
- {|\Delta_1|}^2$ is another invariant.  The significance of this can
be seen as follows.  Pick $\omega\in P$ and let $e=e(\omega)$ be the
corresponding real basis.  Then $d\om3\wedge d\om3 = -8\Delta_1
e^{1234}$ and $d\om3\wedge d\bom3 = -8\Delta_2 e^{1234}$, whence
\begin{gather*}
  de^5\wedge de^5 =
  -4\left(\Delta_2 + \re(\Delta_1) \right)e^{1234},\\
  de^6\wedge de^6 =
  - 4\left(\Delta_2 -\re(\Delta_1) \right)e^{1234},\\
  de^5\wedge de^6 = -4\im(\Delta_1)e^{1234}.
\end{gather*}
The numbers $\Delta_2 \pm \re(\Delta_1)$ determine whether or not the
two-form $de^5$ and $de^6$ are \emph{simple} or not.  A two-form
$\alpha$ is simple if and only if $\alpha\wedge\alpha=0$. The equation
\begin{gather}
  \label{eq:28}
  (s de^5 - t de^6)\wedge(s de^5 - t de^6)=0.
\end{gather}
 is equivalent to the second order homogeneous equation
\[
 (\Delta_2 +\re(\Delta_1))s^2
 - 2\im(\Delta_1) st + (\Delta_2 -\re(\Delta_1)) t^2 = 0.
\]
As the discriminant of this equation is $ {|\Delta_1
|}^2-\Delta_2^2$, it has non-trivial real solutions if and only if
${|\Delta_1 |}^2-\Delta_2^2\ge 0$.

If $d\om3$ and $d\bom3$ are linearly independent, a solution $(s,t)$
to~\eqref{eq:28} exists precisely when $sde^5+tde^6$ is simple.
When ${|\Delta_1 |}^2-\Delta_2^2=0$ there is precisely one such
non-trivial solution, when ${|\Delta_1|}^2-\Delta_2^2>0$ there are
two.  When $d\om3$ and $d\bom3$ are linearly dependent it is easy to
see from equations~\eqref{eq:32}, \eqref{eq:35} and \eqref{eq:36}
that ${|\Delta_1 |}^2=\Delta_2^2$.

\subsection{Identification of underlying real algebras}
Given the invariants of the last section, we now have the means filter isomorphism classes of $\lie g$ for
a given set of structure constants $\epsilon,\rho,A,B,C,D$ of a
nilpotent complex structure. As we determine the underlying real
algebras, we also identify all the invariants in the complex
structure equations in the next few paragraphs.

\begin{lemma}\label{lem:4} The following
  statements are equivalent.
  \begin{enumerate}[$(1)$]
  \item For every nilpotent complex structure $J$ on $\lie g$ and
    every $\omega$ in $P(\lie g,J)$, the condition $\Delta_2(\omega) = 0 = \Delta_1(\omega)$
    holds. \label{item:5}
  \item There exists a nilpotent $J$ on $\lie g$ and some $\omega\in P(\lie
    g,J)$ such that $\Delta_2(\omega) = 0 =
    \Delta_1(\omega)$.\label{item:6}
  \item The Lie algebra $\lie g$ is isomorphic to one of the following
    \begin{gather*}
      \lie h_1 = (0,0,0,0,0,0),\quad \lie h_8=(0,0,0,0,0,12), \quad
      \lie h_6=(0,0,0,0,12,13), \\ \lie h_7=(0,0,0,12,13,23),\quad
      \lie h_{10}=(0,0,0,12,13,14),\quad \lie h_{16}=(0,0,0,12,14,24).
    \end{gather*}\label{item:7} \end{enumerate} \end{lemma}
\bproof It is clear that \eqref{item:5} implies \eqref{item:6}.  Now
suppose \eqref{item:6} holds: pick $J$ and $\omega$ so that
$\Delta_2(\omega) = 0 = \Delta_1(\omega)$.  Since
$d\om2,d\bom2,d\om3,d\bom3$ span $d\lie g^*_c$ and $d\om2\wedge d\om3
= 0 = d\om2\wedge d\bom3$ by the nilpotency of $J$, any two elements
$\alpha_1,\alpha_2$ in $d\lie g^*_\CC$ satisfy
$\alpha_1\wedge\alpha_2=0$.  Since this is in particular also true for
the real elements, a basis of simple two-forms for $d\lie g^*$ exist
so that any two basis elements satisfy $\alpha_1\wedge\alpha_2=0$.
Now consult the classification of six dimensional nilpotent Lie
algebras with complex structures \cite[Theorem 2.9]{Ugarte}.  This
gives \eqref{item:7}.  If \eqref{item:7} holds then any $\omega$ in
$P(\lie g,J)$ for any complex structure $J$ on $\lie g$ has
$d\om{i}\wedge d\om{j}=0=d\om{i}\wedge d\bom{j}$ for all $i,j$.  This
completes the proof.  \eproof

\begin{corollary}\label{cor:2}
  Suppose $\lie g$ is not one of the Lie algebras listed in Lemma{\rm
    ~\ref{lem:4}}.  For any integrable nilpotent $J$ and any $\omega$
  in $P(\lie g,J)$, one has $\Delta_2(\om{} )^2+\abs{\Delta_1(\om{}
    )}^2>0$.
\end{corollary}

\begin{lemma}\label{lem:5}
  The following statements are equivalent.
  \begin{enumerate}[$(1)$]
  \item For every nilpotent complex structure $J$ on $\lie g$ and
    every $\omega$ in $P(\lie g,J)$, the condition $\Delta_2(\omega)^2 <
    \abs{\Delta_1(\omega)}^2$ holds. \label{item:8}
  \item There exists a nilpotent $J$ on $\lie g$ and some $\omega$ in $P(\lie
    g,J)$ such that the inequality $\Delta_2(\omega)^2<
    \abs{\Delta_1(\omega)}^2$ is satisfied.\label{item:9}
  \item The Lie algebra $\lie g$ is isomorphic to one of the following
    \begin{gather*}
      \lie h_2 = (0,0,0,0,12,34), \quad
      \lie h_{12}=(0,0,0,12,13,24), \quad \lie
      h_{13}=(0,0,0,12,13+14,24).
    \end{gather*}\label{item:10}
  \end{enumerate}
\end{lemma}
\bproof The implication \eqref{item:8}$\Rightarrow$\eqref{item:9} is
trivial.  Suppose that $J$ and $\omega$ are given as in
\eqref{item:9}.  Solving \eqref{eq:28}, we get two real, simple
two-forms in the span of $d\om3, d\bom3$.  It follows that $d\lie g^*$
has a basis consisting only of simple two-forms.  The classification
\cite[Theorem 2.9]{Ugarte}, Lemma~\ref{lem:4} and
Corollary~\ref{cor:2} give \eqref{item:10}.

Now suppose that \eqref{item:10} holds and let $J$ be a nilpotent
complex structure on $\lie g$.  Pick any $\omega$ in $P(\lie g,J)$.
Represent $\lie h_2$ as $(0,0,0,0,13,24)$. For any of the three
algebras listed, any nilpotent complex structure and any $\omega$ in
$P(\lie g,J)$, there are constants $a,b,c,r$ such that
$d\om3=ae^{12} + b(e^{13}+re^{14})+ c e^{24}$ where $r=0$ or $1$. So
$d\om3\bw d\om3 = -2bc e^{1234}$ and $d\om3\bw d\bom3 = -(b\bar c
+\bar b c) e^{1234}$.  By Corollary~\ref{cor:2}, $bc\not=0$ so
(after re-scaling of $\om1$ and $\om2$) we have $\abs{\Delta_1}^2 =
\abs{b\bar c}^2 \ge \re(b\bar c)^2 = \tfrac14\abs{b\bar c + \bar b
  c}^2 = \Delta_2^2$.  Equality occurs precisely if $b\bar c$ is real.

To see that this does not occur, note that by nilpotency of $J$,
$d\om2\bw d\om2=0=d\om2\bw d\om3$, whence $d\om2 = ue^{12}$ for some
complex number $u$.  If $u=0$ then $\lie g=\lie h_2$ and $a=r=0$.
Otherwise, take $u=1$ and $\om3-\om2$ as a `new' $\om3$.  This has
$a=0$. So for all three algebras and all $J$, we can take an $\omega$
in $P(\lie g,J)$ with $a=0$.  Then $d\om3$ and $d\bom3$ are linearly
dependent precisely when $b\bar c$ is real.  In this case $n_1$ is $4$
if $\epsilon\not=0$, and $5$ otherwise.  The latter value is not
realized for the given algebras.  Only $(0,0,0,0,13,24)$ has $n_1=4$
but clearly $d\om2\bw d\om2=0=d\om2\bw d\om3$ shows that for this
algebra $\epsilon=0$ for all $J$.  Therefore $b\bar c$ is never real
and so $\abs{\Delta_1}^2>\Delta_2^2$ \eproof

\begin{lemma}\label{lem:6}
  The following
  statements are equivalent.
  \begin{enumerate}[$(1)$]
  \item For every nilpotent complex structure $J$ on $\lie g$ and
    every $\omega\in P(\lie g,J)$, the condition $\Delta_2(\omega)^2 >
    \abs{\Delta_1(\omega)}^2$ holds. \label{item:11}
  \item There exists a nilpotent $J$ on $\lie g$ and some $\omega\in P(\lie
    g,J)$ such that the inequality $\Delta_2(\omega)^2>
    \abs{\Delta_1(\omega)}^2$ is satisfied.\label{item:12}
  \item The Lie algebra $\lie g$ is isomorphic to one of the following
    \begin{gather*}
      \lie h_5 = (0,0,0,0,13+42,14+23), \quad
      \lie h_{15}=(0,0,0,12,13+42,14+23).
    \end{gather*}\label{item:13}
  \end{enumerate}
\end{lemma}
\bproof The idea is as for the preceding Lemmas.
Suppose~\eqref{item:12}.  There are then no simple elements in the
real span of $d\om3+d\bom3,~i(d\om3-d\bom3)$ as
equation~\eqref{eq:28} has no real solutions. This of course means
that for the real basis $e(\omega)$ all linear combinations of
$de^5$ and $de^6$ are non-simple.  This also holds for all
elements in the span of $de^4,de^5,de^6$.  In \cite[Theorem
2.9]{Ugarte} only two algebras have the property that all elements
in the span of $\{de^i\}$ are non-simple.  These are listed in
(\ref{item:13}).

  Building a nilpotent $J$ from $\lie
  h_5$ or $\lie h_{15}$ gives $d\om3=a e^{12} +
  b(e^{13}+e^{42})+c(e^{14}+e^{23})$.  Then $\Delta_1 = b^2+c^2$ and
  $\Delta^2_2 = \abs{b}^2+\abs{c}^2$, so $\Delta_2^2-\abs{\Delta_1}^2 =
  2(\abs{b\bar c}^2-\re{(b\bar c)^2})= 2\im(b\bar c)\ge 0$, with equality if
  and only if $b\bar
  c$ is real.  Arguing as in the proof of Lemma \ref{lem:5} one shows
  that $b\bar c$ cannot be real.  It proves the implications
  \eqref{item:12} to \eqref{item:13}.  The
  implication \eqref{item:11} to \eqref{item:12} is obvious.
\eproof

Now one case is left, namely $\abs{\Delta_1}^2 = \Delta_2^2>0$. By
Lemmas~\ref{lem:4},~\ref{lem:5} and~\ref{lem:6} this condition must
characterize the remaining algebras in the classification of
Lemma~\ref{lem:3}.  In one special case we may be more explicit.
\begin{lemma}\label{lem:7}
  Suppose $\omega$ in  $P(\lie g,J)$ is such that $\abs{\Delta_1}^2 =
\Delta_2^2>0$ and $d=1$.  Then there exists  $\sigma$ in $P(\lie
g,J)$
  such that the real basis $e(\sigma)$ has the following structure
  equations
  \begin{itemize}
  \item $\lie h_3=(0,0,0,0,0,12-\sign(\Delta_2) 34)$ if
    $\epsilon=0,~\abs{\Delta_1}^2=\Delta_2^2>0$,
  \item $\lie h_9=(0,0,0,12,0,14-23)$ if $\epsilon=1$ and
    $\abs{\Delta_1}^2=\Delta_2^2>0$.
  \end{itemize}
\end{lemma}
\bproof Suppose that $\epsilon=0$ and $\abs{\Delta_1}^2=\Delta_2^2>0$.
Note that $\Delta_2=\abs{C}^2-\bar A D$ by~\eqref{eq:32}.  Define
$\lambda>0$ by $\Delta_2 = \sign(\Delta_2)\lambda^2$.  Then
  \begin{equation*}
    \bar A d\om3 = (A\om1 + C\om2)(\bar A\bom1+\bar C\bom2) -
    \sign(\Delta_2)\lambda^2\om2\bom2,
  \end{equation*}
  which gives that second part if $A\not=0$.  If $D\not=0$ we rewrite
  similarly.

  If $A=0=D$, $\Delta_2=\abs{B}^2=\abs{C}^2>0$. We
  pick square roots of $B$ and  $C$, and
  set
  \[
  \te1=\frac{1}{\sqrt2}\left(\sqrt{\frac{B}{C}}\om1+\om2\right), \quad
  \te2=\frac{1}{\sqrt2}\left(-\om1+\sqrt{\frac{C}{B}}\om2\right),
  \quad \te3= -\frac{1}{2\sqrt{BC}}\om3.
  \]
    Then $d\te3 = -(1/2) (\te1\bte1 -
  \te2\bte2)$, whence $de^6=e^{12}-e^{34}$.

  If $\epsilon=1$ then  $A=D=0$. If $\abs{\Delta_1}^2=\Delta_2^2>0$, we take
\[
\te1=2\sqrt{\frac{B}{C}}\om1, \quad \te2=-2\om2, \quad
\te3=\frac{2}{\sqrt{BC}}\om3
\]
  to get $d\te1 = 0, ~d\te2 =
  -(1/2)\te2,~d\te3 = -(1/2)(\te1\bte2 + \te2\bte1)$.
\eproof

\begin{lemma}
  Suppose $\omega$ in $P(\lie g,J)$ is such that $\abs{\Delta_1}^2 =
  \Delta_2^2>0$, $d\om3$ and $d\bom3$ are linearly independent.  Then
  $\lie g$ is one of $\lie h_4, {\lie h}_{11}, {\lie h}_{14}$, with
  $n_2$ being an invariant to distinguish the different spaces.
\end{lemma}
\bproof The proof is similar to the one of the last lemma. As this
is the last remaining case in the classification of all nilpotent
complex structures, one may also identify the concerned algebras
using \cite{Salamon} or \cite{Ugarte}.  \eproof

Next, we tabulate the invariants for all nilpotent complex
structures according to their underlying nilpotent algebras.
\begin{theorem}
  A nilpotent complex structure on a six-dimensional nilpotent Lie
  algebra $\lie g$ is determined by the data of its nilpotent complex
  structures, and vice-versa, as indicated in Table \ref{tab:1} below.
\end{theorem}


\begin{table}[h]
  \begin{center}
    \begin{tabular}{|@{}|l|l|c|c|c|c|c|c|c|@{}|}
      \hline \( n \) & \( \lie g \) & \( \abs{\Delta_1}^2-\Delta_2^2 \) &
      \( \abs{\Delta_1} \) & \( \abs{\Delta_2} \) & \( {\epsilon} \)
      & \( \abs{\rho}  \) & \( d \) \\
      \hline \( (6,6) \) & \( \lie h_1=(0,0,0,0,0,0) \)
      & \( 0 \) & \( 0 \) & \( 0 \) & \( 0 \) & \( 0 \) & \( 0 \) \\
      \( (5,6) \) & \( \lie h_8=(0,0,0,0,0,12) \)
      & \( 0 \) & \( 0 \) & \( 0 \) & \( 0 \) & \( 0 \) & \( 1 \) \\
      \( (5,6) \) & \( \lie h_3=(0,0,0,0,0,12+34) \)
      & \( 0 \) & \( + \) & \( + \) & \( 0 \) & \( 0 \) & \( 1 \) \\
      \( (4,6) \) & \( \lie h_6=(0,0,0,0,12,13) \)
      & \( 0 \) & \( 0 \) & \( 0 \) & \( 0 \) & \( + \) & \( 2 \) \\
      \( (4,6) \) & \( \lie h_4=(0,0,0,0,12,14+23) \)
      & \( 0 \) & \( + \) & \( + \) & \( 0 \) & \( * \) & \( 2 \) \\
      \( (4,6) \) & \( \lie h_2=(0,0,0,0,12,34) \)
      & \( + \) & \( + \) & \( * \) & \( 0 \) & \( * \) & \( 2 \) \\
      \( (4,6) \) & \( \lie h_5=(0,0,0,0,13+42,14+23) \)
      & \( - \) & \( * \) & \( + \) & \( 0 \) & \( * \) & \( 2 \) \\
      \hline \( (4,5) \) & \( \lie h_9=(0,0,0,0,12,14+25) \)
      & \( 0 \) & \( + \) & \( + \) & \( 1 \) & \( 0 \) & \( 1 \) \\
      \( (3,6) \) & \( \lie h_7=(0,0,0,12,13,23) \)
      & \( 0 \) & \( 0 \) & \( 0 \) & \( 1 \) & \( + \) & \( 2 \) \\
      \( (3,5) \) & \( \lie h_{10}=(0,0,0,12,13,14) \)
      & \( 0 \) & \( 0 \) & \( 0 \) & \( 1 \) & \( + \) & \( 2 \) \\
      \( (3,5) \) & \( \lie h_{11}=(0,0,0,12,13,14+23) \)
      & \( 0 \) & \( + \) & \( + \) & \( 1 \) & \( + \) & \( 2 \) \\
      \( (3,5) \) & \( \lie h_{12}=(0,0,0,12,13,24) \)
      & \( + \) & \( + \) & \( * \) & \( 1 \) & \( + \) & \( 2 \) \\
      \( (3,4) \) & \( \lie h_{16}=(0,0,0,12,14,24) \)
      & \( 0 \) & \( 0 \) & \( 0 \) & \( 1 \) & \( + \) & \( 2 \) \\
      \( (3,4) \) & \( \lie h_{13}=(0,0,0,12,13+14,24) \)
      & \( + \) & \( + \) & \( * \) & \( 1 \) & \( + \) & \( 2 \) \\
      \( (3,4) \) & \( \lie h_{14}=(0,0,0,12,14,13+24) \)
      & \( 0 \) & \( + \) & \( + \) & \( 1 \) & \( + \) & \( 2 \) \\
      \( (3,4) \) & \( \lie h_{15}=(0,0,0,12,13+24,14+23) \)
      & \( - \) & \( * \) & \( + \) & \( 1 \) & \( * \) & \( 2 \) \\
      \hline
    \end{tabular}
  \end{center}
  \caption { $\lie g$ and parameters in the complex structure
    equations.  In the table, \lq$\ 0$', \lq$+$' and \lq$-$' indicates that the
    value of the corresponding number is zero, positive or negative,
    while \lq$*$' means that the value is constrained only by the data
    to its left in the table.  The number $d$ of the right-most column
    is the dimension of the linear span of $d\om3$ and $d\bom3$.} \label{tab:1}
\end{table}


\begin{remark}\label{rem:7}\  \rm
  It is known that each of the four algebras with a \lq$*$' in the
  $\abs{\rho}$ column admits both abelian and non-abelian complex
  structures~\cite{Ugarte}.  For $\lie h_5$ and $\lie h_{15}$ this is
  particularly easy to see as both may be represented with either
  $d\om3=\om1\om2$ or $d\om3=\om1\bom2$.  An abelian complex
  structures on $\lie h_2$ is given by $d\om3=i\om1\bom1+\om2\bom2$
  and on $\lie h_4$ by $d\om3=i\om1\om1 + \om1\bom2+\om2\bom2$.  A
  non-abelian nilpotent complex structure on $\lie h_2$ and $\lie h_4$
  may be obtained for instance by setting $d\om1=0=d\om2$ and $d\om3 =
  \rho\om1\om2+ B\om1\bom2+B^{-1}\om2\bom1$ for some $B$ such that
  $\abs{B}\not=1$ with $\abs{\rho}^2=(\abs{B}\pm\abs{B^{-1}})^2$ for
  $\lie h_4$ and
  $(\abs{B}-\abs{B^{-1}})^2<\abs{\rho}^2<(\abs{B}+\abs{B^{-1}})^2$ for
  $\lie h_2$.  We note that any other choice of $\rho$ gives a
  non-abelian complex structure on $\lie h_5$ and one on $\lie h_{15}$
  if we take $d\om2=\om1\bom1$ instead.

  For the algebras with $*$'s in the other columns, i.e. those with
  different values of $\abs{\Delta_1}^2,~\Delta_2^2$, it is also
  always possible to find a complex structure such that the smaller of
  the two is zero.
\end{remark}

\begin{lemma} If  $d=2$,
   $\Delta_1(\omega)=0=\Delta_2(\omega)$, then the complex
  structure is non-abelian.

  Furthermore,
$\epsilon=0$ if and only if there exists $\sigma$ in $ P(\lie g,J)$
such that
  $e(\sigma)$ has structure equations $\lie h_6=(0,0,0,0,13,14)$.
  When $\epsilon=1$, one of the following three cases occurs.
  \begin{itemize}
  \item If there exists an $\omega$ in $P(\lie g,J)$ such that $B=0$,
    there exists a $\sigma$ such that the equation for appropriate
    $e(\sigma)$ is ${\lie h}_{10}=(0,0,0,12,13,14)$.
  \item If there exists an $\omega$ in $P(\lie g,J)$ such that
    $B/\rho>0$, a $\sigma$ may be chosen such that the equation for
    appropriate $e(\sigma)$ is $\lie h_7=(0,0,0,12,13,23)$.
  \item Otherwise $\sigma$ may be chosen such that the equation for
    appropriate $e(\sigma)$ the structure equations are $\lie
    h_{16}\cong(0,0,-t(12),s(12),13,23)$ and $(s+it)^2=B/\rho$.
  \end{itemize}
\end{lemma}
\bproof A simple exercise in algebra using the expressions
\eqref{eq:35} and \eqref{eq:36}, Lemma \ref{lem:3} and Lemma
\ref{dependence}(\ref{item:14}) shows that if
$\Delta_1(\omega)=0=\Delta_2(\omega)$ and $\rho=0$ then $d\om3$ and
$d\bom3$ are linearly dependent. This gives the first statement.

  Suppose that $\epsilon=0$. If in addition $A=0$, then
  $\Delta_1=-BC=0$.  When $B=0$ and
  $\abs{C}^2=\abs{\rho}^2$.  Then we may rearrange to get
  \begin{align*}
    d\om3=(\rho(\om1+(\bar D/\bar C)\om2) - C(\bom1+(D/C)\bom2))\om2.
  \end{align*}
  So choose $r,c$ such that $r^2=\rho$ and $c^2=-C$ and set
  $\te1=(r/c)(\om1+(\bar D/\bar C)\om2)$, $\te2=2\om2$ and
  $\te3=\om3/(cr)$ we get $d\te3=\tfrac12(\te1+\bte1)\te2$.
  If $C=0$, we note that
  \begin{align*}
    d\om3=(\om1+D/B\om2)(\rho\om2 + B\bom2).
  \end{align*}
  Take $r,b$ such that $r^2=\rho$ and $b^2=B$ and set
  $\ta1=-(r/c)\om2,~\ta2 = 2(\om1+D/B\om2),~\ta3 =\om3/(bc)$.  Then
  $d\te3=\tfrac12(\te1+\bte1)\te2$, again.
  If $D=0$ instead of $A=0$, we interchange $\om1$ and $\om2$ and proceed with an argument as above.

  Finally, if $AD=BC\not=0$, we may write
  \[
    d\om3 = ((\abs{C}^2 - \bar A D)(\om1+(\bar D/\bar C)\om2) +
    C(\bom1+(D/C)\bom2))((A/C)\om1+\om2).
  \]
  Since $0<\abs{\rho}^2=\abs{\abs{C}^2- \bar A D}^2/\abs{C}^2$ this is
  equivalent to $d\te3=\tfrac12(\te1+\bte1)\te2$.

  When $\epsilon=1$, then $A=0=D$. As $\Delta_1=0$, by definition
  (\ref{eq:35}) $BC=0$.  If $B=0$, then $d\om3=(\rho\om1-C\bom1)\om2$,
  which we may treat precisely as above to get $d\te1=0,
  ~d\te2=-\tfrac12\te1\bte1, ~d\te3=\tfrac12(\te1+\bte1)\te2$.  If
  $C=0$, pick square roots: $r^2=\rho,~b^2=B$ and set
  $\te1=\om1,~\te2=-\tfrac12(r/b)\om2, ~\te3=\tfrac12\om3$.  Then
  $d\te3 = \te1(\te2+\bte2)$ but $d\te2=(r/b)\te1\bte1$.  Writing
  $r/b= s+it$ we get
  \[
      de^1 =0, \quad de^2=0,\quad
      de^3=-t e^{12}, \quad
      de^4=se^{12},\quad
      de^5= e^{13},\quad
      de^6= e^{23}.
\]
  When $r/b$ is real (which happens if and only if $\rho/B>0$) this is precisely
  $(0,0,0,12,13,23)$.  When $r/b$ is purely imaginary, we get
  $(0,0,12,0,13,23)\cong\lie h_{16}$.  Otherwise, replace $e^4$ by $se^3+te^4$ and
  divide $e^3, e^5$ and $e^6$ with $-t$ to get $(0,0,12,0,13,23)$
  again. 
\eproof

\begin{corollary}\label{cor:22}
  There are no abelian complex structures on $\lie h_p$ for
  $p=6,7$,$10,11$,$12,13$,$14,16$.
  Moreover, suppose that $\omega$ in $P(\lie
  g,J)$ has structure constants $\epsilon,\rho,A,B,C,D$.  If
  $\epsilon=0=\rho$ and $\Delta_1=0$, then $\Delta_2\ge 0$ with $\Delta_2=0$
  if and only if $d\om3$ and $d\bom3$ are linearly dependent.  If
  $\epsilon=1$ and $\rho=0$ then $\Delta_2^2-\abs{\Delta_1}^2\ge 0$
  with equality if and only if $d\om3$ and $d\bom3$ are linearly
  dependent.
\end{corollary}
\bproof
  For $p=7,10,11$ and $12$ this was established by Lemma~\ref{lem:3}.
  For $p=6$ and $16$, any complex structure on $\lie h_p$ has
  $\Delta_1=0=\Delta_2$ by Lemma~\ref{lem:4}.  However, $\Delta_1=0$
  with $\rho=0=\epsilon$ implies $\Delta_2\ge0$ with equality if and
  only if $d\om3$ and $d\bom3$ are linearly dependent.  For $p=13$ and
  $14$, $\epsilon=1$ and $\abs{\Delta_1}^2\ge\Delta_2^2$.  The first
  statement may then be seen to follow from the second and third.

  If $\epsilon = 0 = \rho = \Delta_1 = 0$ then clearly
  \begin{equation*}
    2\Delta_2 =
    \begin{cases}
      \abs{B}^2+\abs{C}^2,&\text{if $A=0$},\\
      \abs{\bar A B- A\bar C}^2/\abs{A}^2,&\text{if $A\not=0$}.
    \end{cases}
  \end{equation*}
  In either case $\Delta_2\ge0$.  If $\Delta_2=0$, $d\om3=D\om2\bom2$
  in the first case, and $\bar A B = A\bar C$ in the second.  It is
  now easy to see that the equations of Lemma \ref{dependence}(\ref{item:14}) are
  satisfied in either case.

  If $\epsilon=1$ and $\rho=0$
  \begin{equation*}
    \Delta_2^2 - \abs{\Delta_1}^2 = (\abs{B}^2-\abs{C}^2)\ge 0
  \end{equation*}
  so equality implies $\abs{B}=\abs{C}$.  Since we may assume that
  $A=0$ when $\epsilon=1$ this shows that $d\om3$ and $d\bom3$ are
  linearly dependent via Lemma \ref{dependence}(\ref{item:14}).  \eproof

\section{Classification of $\DGA ({\lie g}, J)$}\label{sec:f}
\label{sec:isom-class-lie} In this section we calculate the
isomorphism class of the six-dimensional complex Lie algebras $\lie
f^1=\lie f^1(\lie g,J)$ obtained from a nilpotent algebra $\lie g$
equipped with a complex structure $J$. Our aim is to identify the
complex Lie algebra structure of $\lie f^1$ for a given $\lie g$ and
$J$. The result will identify $\lie f^1$ as the complexification of
one of the real nilpotent algebras ${\lie h}_n$.

When a complex structure $J$ is given, recall that the Lie algebra
structure on ${\overline{\lie f}}^1$ is defined by
$\dbar{}\colon\lie{g}^{(1,0)}\oplus\lie g^{*(0,1)}\to \Lambda^2(\lie
g^{(1,0)}\oplus\lie g^{*(0,1)})$. If $X\in\lie g^{(1,0)}, {\overline
  Y}\in\lie g^{(0,1)},\omega\in\lie g^{*(0,1)}$, then $\dbar\omega$ is
the (2,0)-component of $d\omega$ and $(\dbar X)({\overline Y})$ is the
(1,0)-part of the vector $- [X,{\overline Y}]^{1,0}$. Let $T_1, T_2,
T_3$ be dual to $\omega^1, \omega^2, \omega^3$. Given the
equations~\eqref{eq:6'}, the differential $\dbar$ is determined by the
following structure equations.
\begin{gather}
  \begin{cases}
    \dbar\bom1 = 0,\quad \dbar\bom2 = 0,\quad \dbar\bom3 = \rho\bom{12},\\
    \dbar T_1 = \epsilon\bom1 T_2 + (A\bom1 + B\bom2) T_3,\quad
    \dbar T_2 = (C\bom1 + D\bom2) T_3,\quad
    \dbar T_3 = 0.
  \end{cases}
\end{gather}
The Schouten bracket is an extension of the following Lie bracket on
${\lie f}^1$.
\begin{gather}
\begin{cases}
[T_1\bullet T_2]=-\rho T_3, \\
[T_1\bullet\bom2]=-\epsilon\bom1, \quad
[T_1\bullet\bom3]=-{\overline A}\bom1-{\overline C}\bom2, \quad
[T_2\bullet\bom3]=-{\overline B}\bom1-{\overline D}\bom2.
\end{cases}
\end{gather}

In this section, we ignore at first the Lie algebra structure on $\lie
f^1$ and focus on the differential structure $\dbar$ of ${\lie f}^1$
seen as a differential graded algebra.  Inspecting the differential
algebra structure, we identify the Lie algebra structure of $({\lie
  f}^1)^*\cong {\overline {\lie f}}^1$ as the complexification of
$\lie h_n$ for some $n$. Taking complex conjugation, we recover the
Lie algebra structure on $\lie f^1$ as a complexification of the same
$\lie h_n$. The results are presented in Table \ref{tab:f1}.

In the presentation below, the subscript $\CC$ in the identification
${\lie f}^1\cong ({\lie h}_n)_\CC$ is suppressed.

Change basis by setting
\begin{equation}\label{eq:11}
  (\et1,\et2,\et3,\et4,\et5,\et6):=(\bom1,T_3,\bom2, T_2,\bom3, T_1).
\end{equation}
This gives the following structure equations
\begin{equation}\label{structure}
    \dbar\eta_1 =\dbar\eta_2=\dbar \eta_3= 0,\quad \quad
    \dbar\et4 = C\et{12}+D\et{32},\quad
    \dbar\et5 = \rho\et{13},\quad
    \dbar\et6 = \epsilon\et{14} + A\et{12} + B\et{32},
\end{equation}
which clearly define a complex $6$-dimensional nilpotent Lie
algebra.

When the invariants $\epsilon, \rho, \Delta_1$ and $\Delta_2$ are
given, we shall use the complex structure equations (\ref{structure})
to identify the Lie algebra underlying ${\overline{\lie f}}^1$ and
hence $\lie f^1$. On the other hand, we use the invariants and the
classification in Table \ref{tab:1} to identify the originating Lie
algebra $\lie g$.  These are listed in the right most column of Table
\ref{tab:f1}.

\subsection{The cases when $\epsilon=0$.} By Corollary~\ref{cor:1}, $n_2=6$.
Then the potentially non-zero structure
equations are
\begin{equation}
    \dbar\et4 = C\et{12}+D\et{32},\quad
    \dbar\et5 = \rho\et{13},\quad
    \dbar\et6 = A\et{12} + B\et{32}.
\end{equation}
There are six possibilities depending on the rank of
$X:=\begin{spmatrix}
  A&B\\C&D
\end{spmatrix}$ and $\rho$.

\subsubsection{When $\rho=0$.}
\begin{enumerate}[$(1)$]
\item If $\rank X=0$ then $\Delta_1=0$ and $\Delta_2=0$. It follows
that $\lie f_1\cong\lie h_1$ and $\lie g=\lie h_1$.
\item If the rank of $X$ is one then $\Delta_1=0$, $\Delta_2\ge 0$ and
  $\lie f^1\cong\lie h_8$.  By Corollary \ref{cor:22}, $d\om3$ and
  $d\bom3$ are linearly dependent if and only if $\Delta_2=0$.
  Therefore, by Table \ref{tab:1} $\lie g=\lie h_8$ when $\Delta_2=0$,
  and $\lie g=\lie h_5$ when $\Delta_2\neq 0$.
\item If $\rank X=2$ and $\rho=0$ then $\Delta_1\not=0$, $\Delta_2$ is
  unconstrained and $\lie f^1\cong\lie h_6$. By Table \ref{tab:1},
  $\lie g=\lie h_2, \lie h_3, \lie h_4$ or $\lie h_5$.
\end{enumerate}
This case accounts for the first four items in Table \ref{tab:f1}.

 \subsubsection{When $\rho\not=0$.}
\begin{enumerate}[$(1)$]
\item If $\rank X=0$ then $\Delta_1=0$ and $\Delta_2>0$. It follows that $\lie
  f_1\cong\lie h_8$ and $\lie g=\lie h_5$.
\item If $\rank X=1$ then $\Delta_1=0$ and $\Delta_2$ is
  unconstrained.  Then $\lie f^1\cong \lie h_6$. However, when the
  value of $\Delta_2$ varies from zero to non-zero, the algebra $\lie
  g$ changes from $\lie h_5$ to $\lie h_6$.
\item If $\rank X=2$ then $\Delta_1\not=0$, $\Delta_2$ is
  unconstrained and $\lie f^1\cong\lie h_7$. The invariants
  $|\Delta_2 |$ and $|\Delta_1|^2-\Delta_2^2$ help to identify the
  three possibilities $\lie h_2, \lie h_4, \lie h_5$ for the algebra
  $\lie g$.
\end{enumerate}

\subsection{The cases when $\epsilon\neq 0$.} We assume that
$\epsilon=1$, $A=D=0$.  Then the potentially non-zero structure
equations are
\begin{equation}
    \dbar\et4 = C\et{12},\quad
    \dbar\et5 = \rho\et{13},\quad
    \dbar\et6 = \et{14} + B\et{32}.
\end{equation}
\subsubsection{When $\rho=0$.} There are three cases (discarding
$B=0=C$).
\begin{enumerate}[$(1)$]
\item If $C=0$ then $\Delta_1=0$, $\Delta_2>0$. It follows that $\lie f^1\cong\lie h_3$ and $\lie g=\lie h_{15}$.
\item If $B=0$ then $\Delta_1=0$ and $\Delta_2>0$. Then $\lie f^1\cong\lie
  h_{17}=(0,0,0,0,12,15)$ and $\lie g=\lie h_{15}$.
\item If $BC\not=0$ then $\lie f^1\cong\lie h_9$. As
$\Delta_1\not=0$, by Corollary \ref{cor:22},
  $\Delta^2_2-\abs{\Delta_1}^2\ge 0$ with equality if and only if  $d\om3$ and
  $d\bom3$ are linearly dependent. It yields two algebras for $\lie
  g$, namely $\lie h_9$ and $\lie h_{15}$.
\end{enumerate}
\subsubsection{When $\rho\not=0$.} There are four cases for $\lie f^1$:
\begin{enumerate}[$(1)$]
\item If $B=0=C$ then $\lie f^1\cong\lie h_6$. As $\Delta_1=0$, $\Delta_2<0$,
 and $\lie g=\lie h_{15}$.
\item If $C=0, B\neq 0$ then $\lie f^1\cong \lie h_4$.
As $\Delta_1=0$ but $\Delta_2$ is unconstrained, by Table
\ref{tab:1}, $\lie g$ could be one of $\lie h_7, \lie h_{16}$ or
$\lie h_{15}$.
\item If $B=0, C\neq 0$,  then $\lie f^1 \cong \lie h_{10}$.
As $\Delta_1=0$, $\Delta_2$ is unconstrained, we get $\lie g=\lie
h_{10}$ if $\Delta_2=0$. Otherwise, we get $\lie g=\lie h_{15}$.
\item If $BC\not=0$ then $\lie f^1 \cong \lie h_{11}$.
$\Delta_1\not=0$, $\Delta_2$ is unconstrained. An inspection of
Table \ref{tab:1} yields the five different algebras $\lie h_{11},
\lie h_{12}, \lie h_{13}, \lie h_{14}$ and $\lie h_{15}$.
\end{enumerate}

To recap all the computations, we have used the invariants of the
complex structural equations to identity both the underlying real Lie
algebra and the structure of the Lie algebra ${\lie f}^1$.  At the
cost of being repetitive, we recall in the following how the
invariants are defined.
\begin{theorem}\label{invariant theorem} 
  Suppose that $\lie g$ is a real six-dimensional nilpotent algebra
  with a nilpotent complex structure $J$. Then there exists a basis
  $\omega^1, \omega^2, \omega^3$ for $\lie g^{*(1,0)}$ such that
\begin{gather}
  \begin{cases}
    d\om1 = 0, \quad
    d\om2 = \epsilon\om1\bom{1},\\
    d\om3 = \rho\om1\om2 + A\om1\bom1 + B\om1\bom2 +
    C\om2\bom1 + D\om{2}\bom2,
  \end{cases}
\end{gather}
where $\epsilon, \rho\in \{0,1\}$. Moreover, let
\begin{eqnarray*}
&&\triangle_1=AD-BC; \quad  \triangle_2=\frac12[|B|^2+|C|^2-A{\bar
D}-{\bar A}D-|\rho|^2];\\
&& d=\dim_\CC \langle d\om3, d{\bar\omega}^3\rangle, \quad
X=\begin{spmatrix}
  A&B\\C&D
\end{spmatrix}
\end{eqnarray*}
be the invariants associated to the structure equations.  Given a real
algebra $\lie g$ in the right-most column, Table \ref{tab:f1} lists
constraints on the values of the invariants that can be realized by a
complex structure $J$ on $\lie g$, as well as the relevant isomorphism
class of the Lie algebra ${\lie f}^1$ in the left-most column.  A
``\(\, * \)'' indicates an un-constrained invariant.
\end{theorem}

\begin{table}[h]
\begin{center}
\begin{tabular}{|@{}|l|c|c|c|c|c|c|c|c|c|c|@{}|}
\hline \( \lie f^1 \) & \( \epsilon \) & \( \abs{\rho} \) & \( \rank
X \) & \( \abs{\Delta_1} \) & \( \abs{\Delta_2} \) & \(
\abs{\Delta_1}^2-\Delta_2^2 \)
& \( d \) & \( \abs{B} \) & \( \abs{C} \) & \( \lie g \) \\
\hline \( \lie h_1 \) & \( 0 \) & \( 0 \) & \( 0 \) & \( 0 \)
& \( 0 \) & \( 0 \) & \( 0 \) &  \( 0 \) & \( 0 \) &  \( \lie h_1 \)  \\
\( \lie h_8 \) & \( 0 \) & \( 0 \) & \( 1 \) & \( 0 \)
& \( 0 \) & \( 0 \) & \( 1 \) & \( * \) &  \( * \) & \( \lie h_{8} \)\\
\( \lie h_8 \) & \( 0 \) & \( 0 \) & \( 1 \) & \( 0 \)
& \( + \) & \( - \) & \( 2 \) & \( * \) &  \( * \) & \( \lie h_{5} \)\\
\( \lie h_6 \) & \( 0 \) & \( 0 \) & \( 2 \) & \( + \)
& \( * \) & \( - \) & \( 2 \) & \( * \) &  \( * \) & \( \lie h_{2}, {\lie h}_3, {\lie h}_4,
{\lie h}_5 \)\\
\hline \( \lie h_8 \) & \( 0 \) & \( + \) & \( 0 \) & \( 0 \)
& \( + \) & \( - \) & \( 2 \) & \( 0 \) &  \( 0 \) & \( \lie h_{5} \)\\
\( \lie h_6 \) & \( 0 \) & \( + \) & \( 1 \) & \( 0 \)
& \( 0 \) & \( 0 \) & \( 2 \) & \( * \) &  \( * \) & \( \lie h_{6} \)\\
\( \lie h_6 \) & \( 0 \) & \( + \) & \( 1 \) & \( 0 \)
& \( + \) & \( - \) & \( 2 \) & \( * \) &  \( * \) & \( \lie h_{5} \)\\
\( \lie h_7 \) & \( 0 \) & \( + \) & \( 2 \) & \( + \)
& \( + \) & \( + \) & \( 2 \) & \( * \) &  \( * \) & \( \lie h_{2} \)\\
\( \lie h_7 \) & \( 0 \) & \( + \) & \( 2 \) & \( + \)
& \( + \) & \( 0 \) & \( 2 \) & \( * \) &  \( * \) & \( \lie h_{4} \)\\
\( \lie h_7 \) & \( 0 \) & \( + \) & \( 2 \) & \( + \)
& \( * \) & \( - \) & \( 2 \) & \( * \) &  \( * \) & \( \lie h_{5} \)\\
\hline \( \lie h_3
\) & \( + \) & \( 0 \) & \( 1 \)
& \( 0 \)
& \( + \) & \( - \) & \( 2 \) & \( + \) &  \( 0 \) & \( \lie h_{15} \)\\
\( \lie h_{17}
\) & \( + \) & \( 0 \) & \( 1 \) & \(
0 \)
& \( + \) & \( - \) & \( 2 \) & \( 0 \) &  \( + \) & \( \lie h_{15} \)\\
\( \lie h_{9} \) & \( + \) & \( 0 \) & \( 1 \) & \( + \)
& \( + \) & \( 0 \) & \( 1 \) & \( + \) &  \( + \) & \( \lie h_{9} \)\\
\( \lie h_{9} \) & \( + \) & \( 0 \) & \( 2 \) & \( * \)
& \( 0 \) & \( - \) & \( 2 \) & \( + \) &  \( + \) & \( \lie h_{15} \)\\
\hline \( \lie h_{6} \) & \( + \) & \( + \) & \( 0 \) & \( 0 \)
& \( + \) & \( - \) & \( 2 \) & \( 0 \) &  \( 0 \) & \( \lie h_{15} \)\\
\( \lie h_{4} \) & \( + \) & \( + \) & \( 1 \) & \( 0 \) & \( 0 \) &
\( 0 \) & \( 2 \) & \( \abs{\rho} \) &  \( 0 \) & \( \lie
h_{7},~\lie h_{16} \)\\
\( \lie h_{4} \) & \( + \) & \( + \) & \( 1 \) & \( 0 \) & \( + \) &
\( - \) & \( 2 \) & \( + \) &  \( 0 \) & \( \lie
h_{15} \)\\
\( \lie h_{10} \) & \( + \) & \( + \) & \( 1 \) & \( 0 \) & \( 0 \)
& \( 0 \) & \( 2 \) & \( 0 \) &  \( \abs{\rho} \) & \( \lie
h_{10} \)\\
\( \lie h_{10} \) & \( + \) & \( + \) & \( 1 \) & \( 0 \) & \( * \)
& \( - \) & \( 2 \) & \( 0 \) &  \( + \) & \( \lie
h_{15} \)\\
\( \lie h_{11} \) & \( + \) & \( + \) & \( 2 \) & \( + \) & \( * \)
& \( + \) & \( 2 \) & \( + \) &  \( + \) & \( \lie
h_{12},~\lie h_{13} \)\\
\( \lie h_{11} \) & \( + \) & \( + \) & \( 2 \) & \( + \) & \( + \)
& \( 0 \) & \( 2 \) & \( + \) &  \( + \) & \( \lie
h_{11},~\lie h_{14} \)\\
\( \lie h_{11} \) & \( + \) & \( + \) & \( 2 \) & \( + \) & \( + \)
& \( - \) & \( 2 \) & \( + \) &  \( + \) & \( \lie
h_{15} \)\\
\hline
\end{tabular}
\end{center}
\caption{$\lie f^1$ as a function of the parameters in the complex
structure equations.} \label{tab:f1}
\end{table}

Ignoring that the same algebra $\lie f^1$ occurs for  distinct
complex structures or different algebras,  we get

\begin{theorem}\label{finding f}
  Given a six-dimensional nilpotent algebra $\lie g$, the associated
  Lie algebra ${\lie f}^1(\lie g, J)$ for all possible nilpotent
  complex structure $J$ are given in the rows of Table
  \ref{tab:g-f1'}.
\end{theorem}
One observes for instance that for $\lie g=\lie h_{15}$ no less than
seven different isomorphism classes are realized for ${\lie f}^1(\lie
g, J)$ as $J$ runs through the space of complex structures on $\lie
g$.  This is a yet another manifestation of the ``jumping phenomenon''
frequently seen in complex structure deformation theory.

Note that the classification of nilpotent Lie algebras in dimension 6
(see \cite{Magnin,Morosov}) over $\mathbb C$ (or $\mathbb R$) has as a
consequence that structure constant may be taken to always be integers,
and in particular real.  Thus any six-dimensional complex nilpotent
algebra is self-conjugate.  Then Proposition \ref{second technical}
implies that the complex isomorphism of Lie algebras between ${\lie
  f}^1$ and ${\lie h}_n$ generates a C-E compatible linear isomorphism
$O\colon\lie h_h\to\lie h_n^*$ such that $\DGA(\lie g, J)$ and
$\DGA(\lie h_n,O)$ are isomorphic as differential
Gerstenhaber algebras. In other words

\begin{theorem}\label{iso dga}
  Given a six-dimensional nilpotent algebra $\lie g$ with a nilpotent
  complex structure $J$, there exists a differential Gerstenhaber
  algebra $\DGA(\lie h, O)$ quasi-isomorphic to $\DGA(\lie g, J)$ if
  and only if the pair $(\lie g,\lie h)$ is checked in Table
  \ref{tab:g-f1'}.
\end{theorem}

\begin{table}[h]
\begin{center}
\begin{tabular}{||l|@{}|c|c|c|c|c|c|c|c|c|c|@{} |}
\hline \( \lie g\backslash \lie f^1({\lie g}, J) \) & \( \lie h_1 \)
& \( \lie h_3
\) & \( \lie h_{4} \) &  \( \lie h_6
\) & \(
\lie h_7 \) & \( \lie h_8 \) & \( \lie h_{9} \) &\( \lie h_{10} \) & \( \lie h_{11} \) & \(\lie h_{17}
\) \\
\hline
\( \lie h_1 \) & \( \tick \) & & & & & & & & & \\
\hline
\( \lie h_2 \) & & & & \( \tick \) & \( \tick \) & & & & & \\
\hline
\( \lie h_3 \) & & & & \( \tick \) & & & & & & \\
\hline
\( \lie h_4 \) & & & & \( \tick \) & \( \tick \) &             & & & & \\
\hline
\( \lie h_5 \) & & & & \( \tick \) & \( \tick \) & \( \tick \) & & & & \\
\hline
\( \lie h_6 \) & & & & \( \tick \) & & & & & & \\
\hline
\( \lie h_7 \) & & & \( \tick \) & & & & & & & \\
\hline
\( \lie h_8 \) & & & & & & \( \tick \) & & & & \\
\hline
\( \lie h_9 \) & & & & & & & \( \tick \) & & & \\
\hline
\( \lie h_{10} \) & & & & & & & & \( \tick \) & & \\
\hline
\( \lie h_{11} \) & & & & & & & & & \( \tick \) &  \\
\hline
\( \lie h_{12} \) & & & & & & & & & \( \tick \) & \\
\hline
\( \lie h_{13} \) & & & & & & & & & \( \tick \) & \\
\hline
\( \lie h_{14} \) & & & & & & & & & \( \tick \) & \\
\hline
\( \lie h_{15} \) & &  \( \tick \) & \( \tick \) & \( \tick \) & & & \( \tick \) & \( \tick \) & \( \tick \) & \( \tick \) \\
\hline
\( \lie h_{16} \) & & & \( \tick \) & & & & & & & \\
\hline
\end{tabular}
\end{center}
\caption{Isomorphism class of $\lie f^1$ against underlying real
  algebra $\lie g$.}
\label{tab:g-f1'}
\end{table}

The algebra $\lie h_{17}$ appears as a candidate for $\lie f^1$ in the
case $\lie g=\lie h_{15}$.  However $\lie h_{17}$ admits no symplectic
structure.  This demonstrates that the differential Gerstenhaber
algebra $\DGA(\lie h, O)$ does not necessarily arise from a symplectic
structure, as remarked at the end of the proof of Proposition
\ref{second technical}.  The issue of whether $\DGA(\lie h, O)$ is or
not coming from a symplectic structure will be deferred to future
analysis.

\section{Application}\label{sec:appl}
Once we identify the Lie algebra structure for ${\lie f}^1({\lie g},
J)$, we have in effect identified the structure of $\DGA({\lie g},
J)$.  Inspired by the concept of weak mirror symmetry \cite{Mer}, one
could well look for seek nilpotent algebras $\lie h$ with symplectic
structure $\Omega$ whose induced differential Gerstenhaber algebra
$\DGA({\lie h}, \Omega)$ is quasi-isomorphic to $\DGA({\lie g}, J)$.
We shall deal with such a general question in the future. At present,
we take advantage of the results in the preceding sections to address
a more focused question.

Supposing that $(J, \Omega)$ is a pseudo-K\"ahler structure on a
six-dimensional real nilpotent algebra $\lie g$, when will there be a
quasi-isomorphism
\begin{equation}\label{self mirror question}
\DGA({\lie g}, J) \rightleftharpoons \DGA({\lie g}, \Omega)\ ?
\end{equation}
Such pseudo-K\"ahler structures can be interpreted as \it weak
self-mirrors\rm, a manifestation of which - in dimension $4$ - was
studied in \cite{Poon}.

In view of Lemma \ref{quasi}, a quasi-isomorphism is in the present
situation equivalent to an isomorphism on the degree-one level:
\[
({\lie f}^1(\lie g, J), [-\bullet -])\cong ({\lie g}^*_\CC,
[-\bullet -]_\Omega) \cong ({\lie g}_\CC, [-\bullet -]).
\]
Recall that a complex structure can be part of a pseudo-K\"ahler
structure on a nilpotent algebra only if it is a nilpotent complex
structure \cite{CFGU}. In view of Table \ref{tab:f1}, a solution
$({\lie g}, J, \Omega)$ for the question (\ref{self mirror question})
could possibly exist only if $\lie g$ is one of the following:
\begin{equation}
\lie h_{1}, \quad \lie h_{6}, \quad \lie h_{8}, \quad \lie h_{9},
\quad \lie h_{10}, \quad \lie h_{11}.
\end{equation}
Below we extract from Table \ref{tab:f1} the invariants for the
candidate complex structures $J$ for these algebras.

\begin{center}
\begin{tabular}{||l|c|c|c|c|c|c|c|c|c|c||}
\hline \( \lie f^1 \) & \( \epsilon \) & \( \abs{\rho} \) & \( \rank
X \) & \( \abs{\Delta_1} \) & \( \abs{\Delta_2} \) & \(
\abs{\Delta_1}^2-\Delta_2^2 \)
& \( d \) & \( \abs{B} \) & \( \abs{C} \) & \( \lie g \) \\
\hline
\( \lie h_1 \) & \( 0 \) & \( 0 \) & \( 0 \) & \( 0 \)
& \( 0 \) & \( 0 \) & \( 0 \) &  \( 0 \) & \( 0 \) &  \( \lie h_1 \)  \\
\hline
 \( \lie h_6 \) & \( 0 \) & \( + \) & \( 1 \) & \( 0 \)
& \( 0 \) & \( 0 \) & \( 2 \) & \( * \) &  \( * \) & \( \lie h_{6}
\)\\
\hline
 \( \lie h_8 \) & \( 0 \) & \( 0 \) & \( 1 \) & \( 0 \)
& \( 0 \) & \( 0 \) & \( 1 \) & \( * \) &  \( * \) & \( \lie h_{8} \)\\
\hline
 \( \lie h_{9} \) & \( + \) & \( 0 \) & \( 1 \) & \( + \)
& \( + \) & \( 0 \) & \( 1 \) & \( + \) &  \( + \) & \( \lie h_{9} \)\\
\hline
 \( \lie h_{10} \) & \( + \) & \( + \) & \( 1 \) & \( 0 \) &
\( 0 \) & \( 0 \) & \( 2 \) & \( 0 \) &  \( \abs{\rho} \) & \( \lie
h_{10} \)\\
 \hline
 \( \lie h_{11} \) & \( + \) & \( + \) & \( 2 \) & \( + \) &
\( + \) & \( 0 \) & \( 2 \) & \( + \) &  \( + \) & \( \lie
h_{11} \)\\
\hline
\end{tabular}
\end{center}

In the next few sections, we shall take the above complex structures,
and seek symplectic structures that realize the quasi-isomorphism
(\ref{self mirror question}). We shall analyze pseudo-K\"ahler
structures on $\lie h_6$, $\lie h_8$, and $\lie h_{11}$ in details,
merely outline the discussion for $\lie h_{9}$ and $\lie h_{10}$, and
skip the trivial case ${\lie h}_1$ completely.

\subsection{${\lie h}_6$}
Given the invariants, the reduced structure equations (\ref{eq:6})
are
\begin{equation}
    d\om1 = 0, \quad
    d\om2 = 0, \quad
    d\om3 = \om1\om2 + A\om1\bom1 + B\om1\bom2 +
    C\om2\bom1 + D\om{2}\bom2.
  \end{equation}
Since $\Delta_1=0$ there exists a constant $\lambda$ such that
either
\[
 d\om3 = \om1\om2 +(\om1+\lambda\om2)(A\bom1+B\bom2) \quad
 \mbox{ or }
 \quad
d\om3 = \om1\om2 +(\lambda\om1+\om2)(C\bom1+D\bom2).
\]
The condition $\Delta_2=0$ implies that in either case, there exists a
change of complex basis so that the structure equations transform to
\begin{equation}\label{standard 6}
  d\om1 = 0, \quad
    d\om2 = 0, \quad
    d\om3 = \om1\om2 + \om1\bom2.
    \end{equation}
It follows that the structure equations for $({\lie f}^1, [-\bullet
-], \dbar)$ are
\begin{equation}
[T_1, T_2]=-T_3, \quad  [T_2, \bom3]=-\bom1, \quad \dbar T_1=\bom2
\wedge T_3, \quad \dbar\bom3=\bom1\wedge\bom2.
\end{equation}

Due to \cite[Lemma 3.4]{CFU}, given the complex structure equations,
any $(1,1)$-form of a compatible symplectic structure is given by
\[
\Omega=a_1\om1\bom1+b_2\om2\bom2+{\overline
a}_2\bom1\om2+a_2\om1\bom2+a_3(\om1\bom3+\bom1\om3),
\]
where $a_1$ and $b_2$ are imaginary numbers and $a_3$ is a real
number. This $2$-form is non-degenerate if and only if $b_2\neq 0$ and
$a_3\neq 0$.

Setting $\om1=e^2+ie^3$, $\om2=-\tfrac12(e^1+ie^4)$ and
$\om3=e^5+ie^6$, reduces the complex structure equation to
\begin{equation}
de^5=e^{12}, \quad de^6=e^{13}.
\end{equation}
Set $a_1=\tfrac{i}2a$, $b_2=2ib$, $a_2=c+ik$ and $a_3=\ell/2$ with
$b\neq 0$ and $\ell\neq 0$. Then the symplectic structure is
\[
\Omega=ae^{23}+be^{14}+c(e^{12}-e^{34})-k(e^{13}+e^{24})+\ell(e^{25}+e^{36}).
\]
Using the contraction with $\Omega$ as an isomorphism from $\lie h_6$ and
$\lie h_6^*$, we obtain a Lie bracket on $\lie h_6^*$ such that
\begin{equation}
  b[e^4, e^5]_\Omega=e^2, \quad b[e^4, e^6]_\Omega=e^3, \quad
  bl[e^5,e^6]_\Omega=(ce^2-ke^3).  
\end{equation}
It is now apparent that the linear map
\begin{equation}
  T_1\mapsto e^5 + \frac{k}{\ell}e^4, \quad T_2\mapsto be^4, \quad T_3\mapsto e^2, \quad
  \bom1\mapsto-e^3, \quad \bom2\mapsto e^1, \quad \bom3\mapsto e^6 + \frac{c}{\ell}e^4.
\end{equation}
yields an isomorphism of differential Gerstenhaber algebras.

Note that the isomorphism exists so long as the symplectic form
$\Omega$ and the designated complex structure $J$ together form a
pseudo-K\"ahler structure.

\begin{proposition}\label{prop h6}
  Let $J$ be any integrable complex structure on $\lie h_6$. Let
  $\Omega$ be any symplectic form on $\lie h_6$ of type $(1,1)$ with
  respect to $J$.  Then the differential Gerstenhaber algebras
  $\DGA(\lie h_6, J)$ and $\DGA(\lie h_6,\Omega)$ are isomorphic.
\end{proposition}

\subsection{${\lie h}_8$}
In this case, the invariants yield the following structure equations.
\begin{equation}
    d\om1 = 0, \quad
    d\om2 = 0, \quad
    d\om3 =A\om1\bom1 + B\om1\bom2 +
    C\om2\bom1 + D\om{2}\bom2,
\end{equation}
where the arrays $(A, B)$ and $(C,D)$ are linearly dependent but are
not identically zero. After a change of complex coordinates, they
could be reduced to
\begin{equation}\label{standard 8}
   d\om1 = 0, \quad
    d\om2 = 0, \quad
    d\om3 =\om1\bom1.
\end{equation}
The induced structure equations for $\lie f^1$ are
\begin{equation}\label{f1-8}
[T_1, \bom3]=-\bom1, \quad \dbar T_1=\bom1\wedge T_3.
\end{equation}
By choosing
\begin{equation}
\om1=e^1+ie^2, \quad \om2=e^3+ie^4, \quad \om3=-2(e^5+ie^6),
\end{equation}
then the real structure equation is indeed the standard one for
$\lie h_8$:
\begin{equation}\label{d6}
de^6=e^{1}\wedge e^2.
\end{equation}

Again, due to \cite[Lemma 3.4]{CFU} given the complex structure
equations, any symplectic $(1,1)$-form is given by
\[
\Omega=ae^{12}+be^{34}+x(e^{13}+e^{24})-y(e^{23}-e^{14})-u(e^{15}+e^{26})+v(e^{25}-e^{16}),
\]
where $a, b, x,y, u,v$ are real numbers. $\Omega$ is non-degenerate
when $b\neq 0$ and $u^2+v^2\neq 0$. Then the induced Lie bracket on
$\lie h_{8}^*$ is given by
\[
[-ue^5-ve^6\bullet ve^5-ue^6]_\Omega=-(ue^2+ve^1). \] Since
$u^2+v^2\neq 0$, it is an elementary exercise to find isomorphism
from $\DGA({\lie h}_8, J)$ to $\DGA({\lie h}_8, \Omega)$. For
instance, when $v\neq 0$, one could construct
 an isomorphism so
that
\begin{equation}
T_1\mapsto -ue^5-ve^6, \quad \bom3\mapsto ve^5-ue^6, \quad
\bom1\mapsto ue^2+ve^1.
\end{equation}

As in the last section, the computation demonstrates more than simply
the existence of a self-mirror pair of complex and symplectic
structure.

\begin{proposition}\label{prop h8} Let $J$ be any integrable complex structure on $\lie h_8$.
Let $\Omega$ be any symplectic form on $\lie h_8$ of type $(1,1)$
with respect to $J$. Then the differential Gerstenhaber algebras
$\DGA(\lie h_8, J)$ and $\DGA(\lie h_8,\Omega)$ are isomorphic.
\end{proposition}

\subsection{${\lie h}_9$}
The complex structure equations are given by
\[
    d\om1 = 0,\quad
    d\om2 = \om1\bom{1},\quad
    d\om3 =   B\om1\bom2 +
    C\om2\bom1,
\]
where $B\neq 0$ and $C\neq 0$. Therefore, we can normalize to
\begin{equation}
  d\om1 = 0,\quad
  d\om2 = -\frac12\om1\bom{1},\quad
  d\om3 =   \frac12\om1\bom2 +\frac12
  \om2\bom1,
\end{equation}
Choose
\begin{equation}
  \om1=e^1+ie^2, \quad \om2=e^4+ie^5, \quad \om3=e^6+ie^3,
\end{equation}
to the effect that $\Omega=e^{13}-e^{26}-e^{45}$ is a pseudo-K\"ahler
form so that ${\lie f}^1(\lie h_9, J)$ is isomorphic to $(\lie h_9^*,
[-\bullet -]_\Omega)$.

\subsection{${\lie h}_{10}$}
The complex structure equation is given by
\begin{equation}
  d\om1 = 0,\quad
  d\om2 = \om1\bom{1},\quad
  d\om3 = \om1\om2  +
  \om2\bom1,
\end{equation}
In this case, when we choose
\[
\om1=e^1+ie^2, \quad \om2=e^3+ie^4, \quad \om3=e^5+ie^6,
\]
then $ \Omega=i(e^{16}-e^{25}-e^{34}) $ is a pseudo-K\"ahler form
such that
    ${\lie f}^1(\lie h_{10}, J)$ is isomorphic to $(\lie h_{10}^*,
    [-\bullet -]_\Omega)$.

\subsection{${\lie h}_{11}$}
This case requires a careful analysis.  We show that \emph{for every
  pseudo-K\"ahler pair $(J,\Omega)$ on $\lie h_{11}$ the differential
  Gerstenhaber algebras $\DGA(\lie h_{11},J)$ and $\DGA(\lie
  h_{11},\Omega)$ are non-isomorphic}.  To this end we shall suppose
that $\Phi\colon\DGA(\lie h_{11},J) \to \DGA(\lie h_{11},\Omega)$ is a
quasi-isomorphism of differential Gerstenhaber algebras obtained from
a pseudo-K\"ahler pair $(J,\Omega)$ and establish a contradiction.

Note that $\lie h_{11}$ is distinguished by the data: $n=(3,5)$ and
$\abs{\Delta_1}^2=\Delta_2^2>0$ for any $J$, see Lemma~\ref{lem:3} and
Lemma~\ref{lem:7}.  Furthermore, for any complex structure on $\lie
h_{11}$ we may always choose a basis of $(1,0)$-forms such that
\begin{equation}
  \label{eq:3}
    d\om1 = 0,\quad
    d\om2 = \om1\bom{1},\quad
    d\om3 =\om1\om2  + B\om1\bom2 +
    C\om2\bom1.
\end{equation}
Choosing $\omega$ this way, the constraints $n_2=5$ and
$\abs{\Delta_1}^2=\Delta_2^2>0$ on the invariants are equivalent to
$B$ being real, $\abs{C}^2=(B-1)^2$ and $BC\not=0$. We shall use this
extensively below.  Precisely these conditions on $B$ and $C$ give
\begin{equation}
  \label{eq:2}
  d((B-1)\om3+C\bom3)=((B-1)\om1+C\bom1)(\om2+\bom2),
\end{equation}
whence 
\begin{equation}
  V_1(\lie h_{11}) = \langle{\om1, \bom1,\om2 + \bom2}\rangle,\qquad
  V_2(\lie h_{11}) = \langle{\om1, \bom1, \om2, \bom2, (B-1)\om3 +
  C\bom3}\rangle.
\end{equation}
Solving the equations $d\Omega=0$ and $\Omega=\bar\Omega$ in the space
of $(1,1)$-forms gives
\begin{equation}
  \label{eq:1}
  \Omega = a_1\om1\bom{1} + a_3(B+1)\om2\bom{2} + a_2\om1\bom2 - \bar
  a_2\om2\bom1 + a_3(\om1\bom3+\om3\bom1),
\end{equation}
where $a_1+\bar a_1=0=a_3+\bar a_3$ and $a_1a_3(B+1)\not=0$ if and
only if $\Omega$ is non-degenerate\footnote{This also means: nilpotent
  complex structures on $\lie h_{11}$ with $B=-1$ have no compatible
  symplectic forms.}.  Therefore $\Omega(T_1) =
a_1\bom1+a_2\bom2+a_3\bom3,~\Omega(T_2) = -\bar a_2\bom1
+a_3(B+1)\bom2,~\Omega(T_3) = a_3\bom1$ and
\begin{gather*}
  \om1=-\frac1{a_3}\Omega(\bar T_3),\qquad \om2 = -
  \frac1{(B+1)a_3}\Omega\left(\bar T_2-\frac{a_2}{a_3}\bar
    T_3\right),\\
  \om3 = -\frac1{a_3}\Omega\left(\bar T_1 + \frac{\bar
      a_2}{(B+1)a_3}\bar T_2 -
    \frac{(B+1)a_1a_3+\abs{a_2}^2}{(B+1)a_3^2}\bar T_3\right).
\end{gather*}
Now the brackets are easily computed
\begin{gather*}
  [\om2\bullet\om3]_\Omega = - \frac1{(B+1)a_3}\om1,\qquad
  [\om2\bullet\bom3]_\Omega
  = - \frac1{(B+1)a_3}\left(\bar C\om1 + B\bom1\right),\\
  [\om3\bullet\bom3]_\Omega = - \frac1{(B+1)a_3^2}\left((a_2+\bar a_2 \bar
    C)\om1 - (\bar a_2+a_2 C)\bom1\right) -
  \frac{B+1}{a_3}(\om2+\bom2),
\end{gather*}
and the lower central series for $(\lie
h_{11}^*,[\cdot\bullet\cdot]_\Omega)$ is
\begin{gather*}
  (\lie h_{11}^*)_1 = \langle \om1, \bom1, \om2 + \bom2\rangle,\quad
  (\lie h_{11}^*)_2 = \langle (B-1)\om1 + C\bom1\rangle,\quad (\lie h_{11}^*)_3=\{0\},
\end{gather*}
while the ascending series is
\begin{gather*}
  D^1(\lie h_{11}^*) = \langle \om1,\bom1\rangle,\quad
  D^2(\lie h_{11}^*) = \langle
  \om1,\bom1,\om2,\bom2\rangle,\quad D^3(\lie h_{11}^*) =
  \lie h^*_{11}.
\end{gather*}
On the other hand, the structure equations for $\DGA(\lie h_{11},J)$
given by~(\ref{eq:3}) are 
\begin{gather*}
  \dbar T_1=\bom1\wedge T_2+B\bom2\wedge T_3, \quad \dbar
  T_2=C\bom1\wedge T_3, \quad \dbar\bom3=\bom1\wedge\bom2,\\
  [T_1\bullet T_2]=-T_3,\quad[T_1\bullet\bom2]=-\bom1,\quad[T_1\bullet\bom3]=-\bar C\bom2,\quad[T_2\bullet\bom3]=-B\bom1.
\end{gather*}
Writing $\lie f^1$ for the space of degree one elements in $\DGA(\lie
h_{11},J)$ we have
\begin{gather*}
  V_1(\lie f^1)=\langle T_3,\bom1,\bom2\rangle,\qquad V_2(\lie f^1) =
  \langle T_2,T_3,\bom1,\bom2,\bom3 \rangle,\\
  (\lie f^1)_1 = \langle T_3,\bom1,\bom2\rangle,\qquad (\lie f^1)_2 =
  \langle \bom1\rangle,\qquad (\lie f^1)_3 = \{0\},\\
  D^1(\lie f^1) = \langle \bom1,T_3\rangle,\quad D^2(\lie f^1) =
  \langle \bom1,T_3,\bom2,T_2\rangle,\quad D^3(\lie f^1) = \lie f^1.
\end{gather*}

By Proposition~\ref{key technical}, any quasi-isomorphism
$\Phi\colon\DGA(\lie h_{11},J) \to \DGA(\lie h_{11},\Omega)$ must be
an isomorphism of DGAs and therefore maps $V_k(\lie f^1)$
isomorphically onto $V_k(\lie h_{11})$, $(\lie f^1)_j$ isomorphically
onto $(\lie h_{11}^*)_j$ and similarly for the ascending sequences.
It follows that complex constants $\phi^{m}_{n}$ exist such that
\begin{eqnarray*}
  \Phi(\bom1)&=&\phi^1_1((B-1)\om1 + C\bom1),\label{phi bom1}\\
  \Phi(T_3) &=& \phi^2_1\om1 +\phi^2_2\bom1, \label{phi t3}\\
  \Phi(\bom2) &=& \phi^3_1\om1 +\phi^3_2\bom1 + \phi^3_3(\om2+\bom2), \label{phi bom2}\\
  \Phi(T_2) &=& \phi^4_1\om1 +\phi^4_2\bom1 + \phi^4_3\om2 +
  \phi^4_4\bom2, \label{phi t2}\\
  \Phi(\bom3) &=& \phi^5_1\om1 +\phi^5_2\bom1 + \phi^5_3\om2 +
  \phi^5_4\bom2 + \phi^5_5((B-1)\om3 + C\bom3), \label{phi bom3}\\
  \Phi(T_1) &=& \phi^6_1\om1 +\phi^6_2\bom1 + \phi^6_3\om2 +
  \phi^6_4\bom2 + \phi^6_5\om3 + \phi^6_6\bom3 \label{phi t1}.
\end{eqnarray*}
More detailed information is now obtained by applying $\Phi$ to the
structure equations.  From $d(\Phi(T_2))=C\Phi(\bom1)\bw\Phi(T_3)$ we
get
\begin{equation}
  \label{eq:8}
  (\phi^4_3-\phi^4_4)=C\phi^1_1((B-1)\phi^2_2-C\phi^2_1).  
\end{equation}
The $((B-1)\om1+C\bom1)(\om2+\bom2)$-component of $d(\Phi(\bom3)) =
\Phi(\bom1)\bw\Phi(\bom3)$ gives
\begin{equation}
  \label{eq:9}
  \phi^5_5=\phi^1_1\phi^3_3.  
\end{equation}
Eliminating $\phi^6_5$ and $\phi^6_6$ in the equations derived from
$d(\Phi(T_1)) = \Phi(\bom1)\bw\Phi(T_2) + B\Phi(\bom2)\bw\Phi(T_3)$
leads to $C\phi^1_1(\phi^4_3-\phi^4_4) +
\phi^3_3((B-1)\phi^2_2-C\phi^2_1)=0$.  The result of inserting
\eqref{eq:8} in this is
$((C\phi^1_1)^2+\phi^3_3)((B-1)\phi^2_2-C\phi^2_1)=0$.  Since $\Phi$ is a
linear isomorphism $\Phi(\bom1)$ and $\Phi(T_3)$ are linearly
independent, and so
\begin{equation}
  \label{eq:10}
  \phi^3_3=-(C\phi^1_1)^2.
\end{equation}
The equation $ [\Phi(T_1)\bullet\Phi(\bom2)]_\Omega=-\Phi(\bom1)$ is
equivalent to 
\begin{equation}
  \label{eq:4}  
  C(B+1)a_3\phi^1_1 = \phi^3_3(C\phi^6_5 - (B-1)\phi^6_6)
\end{equation}
while the $\om2+\bom2$-component of $[\Phi(T_1) \bullet
\Phi(\bom3)]_\Omega = -\bar C\Phi(\bom2)$ gives
\begin{equation}
  \label{eq:5}
  a_3\bar C\phi^3_3 =
(B+1)\phi^5_5(C\phi^6_5 - (B-1)\phi^6_6).
\end{equation}
Substituting first \eqref{eq:9}, and then equations \eqref{eq:4} and
\eqref{eq:10} in \eqref{eq:5} yields
\begin{equation*}
  a_3\abs{C}^2C (\phi^1_1)^2 = - a_3(B+1)^2C(\phi^1_1)^2.
\end{equation*}
Since $\abs{C}^2=(B-1)^2$, this implies $a_3C\phi^1_1=0$ and so
establishes our contradiction: if $a_3=0$ then $\Omega$ is degenerate,
$C=0$ cannot be realized on $\lie h_{11}$, and if $\phi^1_1=0$ then
$\Phi$ is no isomorphism. \eproof

\subsection{Conclusion}

The computation in the past few paragraphs is summarized in the
following observation.

\begin{theorem}\label{main} A six-dimensional nilpotent algebra $\lie g$ admits a
pseudo-K\"ahler structure $(J, \Omega)$ such that $\DGA({\lie g},
J)$ is quasi-isomorphic to $\DGA({\lie g}, \Omega)$ if and only if
$\lie g$ is one of $\lie h_1$, $\lie h_6$, $\lie h_8$, $\lie h_9$
and $\lie h_{10}$.
\end{theorem}
\begin{remark}
  \rm In this paper, we have dealt exclusively with Lie algebras.
  However, it is possible to extend the whole discussion to
  nilmanifolds $M=G/\Gamma$, i.e.  quotients of simply connected
  nilpotent Lie groups $G$ with respect to co-compact lattices
  $\Gamma$. Indeed, the de Rham cohomology of $M$ is given by
  invariant forms on $G$ \cite{Nomizu}. Therefore, when $M$
  has an invariant symplectic structure, the invariant differential
  Gerstenhaber algebra $\DGA(\lie g, \Omega)$ provides a minimal model
  for the differential Gerstenhaber algebra over the space of sections
  of the exterior differential forms on the nilmanifold $M$.

  Similarly, for nilpotent complex structures on nilmanifolds there
  are partial results proving that the space of invariant sections is
  a minimal model of the Dolbeault cohomology with coefficients in the
  holomorphic tangent sheaf \cite{CF} \cite{CFP} \cite{CFGU}.  Given
  such a result for particular class of complex structures (e.g.
  abelian complex structure \cite{CFP}), Theorem \ref{main} can be
  paraphrased as a statement about quasi-isomorphisms of DGAs over
  nilmanifolds with pseudo-K\"ahler structures.
\end{remark}

\begin{ack}
  We are grateful to S. Chiossi for reading the manuscript and for his
  extremely useful comments.
\end{ack}

\end{document}